\documentclass[journal]{IEEEtran}
\usepackage[margin=0.7in]{geometry}

\usepackage{graphicx}      

\usepackage{amsfonts,latexsym,amssymb,amsmath,url,stackengine, stmaryrd}
\usepackage[colorlinks=true, allcolors=blue]{hyperref}
\usepackage{color}
\usepackage{caption}
\usepackage{subcaption}

\usepackage[french,english]{babel} 
\usepackage{csquotes}
\usepackage[T1]{fontenc} 

\catcode`\@=11
\def\downparenfill{$\m@th\braceld\leaders\vrule\hfill\bracerd$}
\def\downparenfill{$\m@th\braceld\leaders\vrule\hfill\bracerd$}
\def\overparen#1{\mskip 2mu\mathop{\vbox{\ialign{##\crcr\crcr \noalign{\kern0.4ex}
\downparenfill\crcr\noalign{\kern0.4ex\nointerlineskip}
$\hfil\displaystyle{#1}\hfil$\crcr}}}\limits\mskip 2mu} 
\catcode`\@=12

\newtheorem{lemma}{Lemma}
\newtheorem{theorem}{Theorem}

\newtheorem{corollary}{Corollary}

\newtheorem{proposition}{Proposition}

\newcommand{\NN}{\mathbb{N}}

\newcommand{\RR}{\mathbb{R}}

\def\T{\mathcal{T}}

\def\C{\mathcal{C}}

\def\L{\mathcal{L}}

\newcommand{\comment}[1]{}

\allowdisplaybreaks

\title{\bf 
Ensembles of Hyperbolic PDEs:\\ Stabilization by Backstepping} 

\author{Valentin Alleaume and Miroslav Krstic

\thanks{V.Alleaume is with the Ecole Mines Paris, PSL Research University, 75006 Paris, France (e-mail: valentin.alleaume@etu.minesparis.psl.eu).}
\thanks{M. Krstic is with the Department
of Mechanical and Aerospace Engineering, University of California at San Diego, La Jolla,
CA, 92093-0411 USA (e-mail: krstic@ucsd.edu).}}

\begin{document}
\maketitle

\begin{abstract}
For the quite extensively developed PDE backstepping methodology for coupled linear hyperbolic PDEs, we provide a generalization 
from finite collections of such PDEs, whose states at each location in space are vector-valued, to previously unstudied infinite (continuum) ensembles of such hyperbolic PDEs, whose states are function-valued. The motivation for studying such systems comes from traffic applications (where driver and vehicle characteristics are continuously parametrized), fluid and structural applications, and future applications in population dynamics, including epidemiology. 
Our design is of an exponentially stabilizing scalar-valued control law for a PDE system in two independent dimensions, one spatial dimension and one ensemble dimension. 
In the process of generalizing PDE backstepping from finite to infinite collections of PDE systems, we generalize the results for PDE backstepping kernels to  the continuously parametrized Goursat-form PDEs that govern such continuously parametrized kernels. The theory is illustrated with a simulation example, which is selected so that the kernels are explicitly solvable, to lend clarity and interpretability to the simulation results.
\end{abstract}

\section{Introduction}

This paper extends the results on control of finitely-many coupled hyperbolic PDEs, using the PDE backstepping method, to control of infinitely many (uncountably many, i.e., a continuum of) such PDEs. 

\paragraph{Backstepping for coupled hyperbolic PDEs} 
Even though PDE backstepping was first developed for parabolic systems \cite{1369395}, the hyperbolic PDE class is nowadays a more natural PDE class from which to begin the study of PDE backstepping. This is because in the easier, hyperbolic case \cite{krstic2008Backstepping,BERNARD20142692},  only a single derivative in space arises. 

The relative ease of dealing with PDEs with a single derivative in space has resulted in the control of hyperbolic PDEs growing into a rich area of research. The greatest progress in the hyperbolic case has been in the development of feedback laws that stabilize coupled hyperbolic systems with fewer inputs than PDEs. A pair of coupled hyperbolic PDEs was stabilized first, with a single boundary input in \cite{Coron2013Local}. An extension to $n+1$ hyperbolic PDEs with a single input was introduced in \cite{dimeglio2013}, an extension to $n+m$ PDEs with boundary actuation on $m$ “homodirectional” PDEs in \cite{hu2016,hu2019boundary}, an extension to cascades with ODEs in \cite{DIMEGLIO2018281},  an extension to “sandwiched” ODE-PDE-ODE systems in \cite{WANG2020109131, 9319184}, and an extension to coupled hyperbolic PDEs with a zero characteristic speed \cite{deandrade2022backstepping}. Redesigns robust to delays are provided in \cite{Auriol2018Delay,Auriol2018Delay1}. PDE backstepping-based output-feedback regulation with disturbances is proposed in \cite{deutscher2018, DEUTSCHER201556}. Adaptive versions of backstepping designs for coupled hyperbolic PDEs are provided in the book \cite{Anfinsen2019Adaptive}. 

We focus here on the extensions of the watershed result in \cite{dimeglio2013}, which considered a finite ensemble of coupled first-order hyperbolic PDEs, to infinite/continuum ensembles of PDEs, under a single input. The capability provided in \cite{dimeglio2013} had had an impact in multiple application domains: for fluid flows, in canals with sediments \cite{diagne2017backstepping} and in canals with unmixed fluids \cite{diagne2017control}, inflexible structures modeled by the Timoshenko model of elasticity  \cite{chen2022backstepping} and in multi-layer beams \cite{chen2023backstepping}, in traffic flows \cite{Yu2019Traffic,burkhardt2021stop,Yu2022} with multiple classes of vehicles (small, medium, trucks) and drivers (young and old, fresh and tired, aggressive and defensive, etc.), and for multi-phase flow in oil drilling (oil, water, natural gas, and “mud” for evacuating rock cuttings during drilling). 

\paragraph{From coupled PDEs, and ensemble ODEs, to ensemble PDEs}
There  is  no  reason  to  ”quantize”  populations  into  unnatural  and  arbitrary  finite  clusters or quantiles when they actually exist on a continuum. Particularly obvious is that vehicles in traffic exist on a continuum relative to weight/inertia and power, and that drivers certainly are on a continuum relative to age or fatigue, drowsiness, and mood. The case for advancing the control of hyperbolic PDE ensembles from finite to continuum/infinite ensembles is further amplified by the fact that, in the continuum limit, the theory becomes notationally cleaner and conceptually clearer. 

We, therefore, consider coupled hyperbolic PDEs on a one-dimensional domain, normalized for clarity to the unit interval $x\in[0,1]$. The ensemble variable $y$ parametrizes the PDEs within the ensemble. For notational convenience we also normalize $y$ so that it belongs to $[0,1]$. In conclusion, we have a PDE system with independent variables$(x,y)\in[0,1]^2$. This may appear as a PDE on a square (2D) domain but it is not. Partial derivatives in space appear only with respect to $x$, with no derivatives in $y$. Only interconnections among the ensemble members appear, as integrals in $y$. Hence, the problem considered is an infinite collection of 1D PDEs, with a scalar/single input at one boundary. Such PDEs are sometimes referred to as ``quasi-2D'' or, in a colloquial discourse, even as ``1.5-D'' PDEs.

From the stabilization perspective, there is a major distinction between 2D and quasi-2D (namely, ensemble 1D) PDEs. While 2D PDEs are, in general, not stabilizable using scalar-valued boundary control (requiring, instead, boundary control that varies not only with time but also with the spatial location on the boundary of the 2D domain), the quasi-2D PDE class we consider lends itself to stabilization using a single PDE input. This drastic difference, primarily from the applicability point of view, is one of the strong motivators for the study of ensemble PDE stabilization.  

One can regard many of the existing results on PDE backstepping control as results for {\em ensemble PDE} systems. The majority of such results are for parabolic PDEs. For instance, the result for coupled parabolic PDEs \cite{vazquez2016boundary}, a distant parabolic cousin of the hyperbolic results in \cite{dimeglio2013,hu2016,hu2019boundary}, is for a finite ensemble of parabolic PDEs. PDE backstepping results in any dimension higher than one, in particular, in 2-D, can be regarded as a very basic form of ensemble PDEs. For instance, the 2D Navier-Stokes PDE governing the infinite-channel (Poiseuille) flow in \cite{vazquez2007closed} is governed by the boundary-controlled 1-D ensemble PDE \cite[eq. (67)]{vazquez2007closed} in the $y$-direction, which we partly quote here for convenience,
\begin{equation}
    u_t = \frac{1}{Re} (u_{yy}-4\pi^2 k^2 u) + 8\pi k i y (y-1) u + \cdots\,,
\end{equation}
and where $k\in \mathbb{R}$, usually referred to as the wavenumber, is the Fourier transform variable in the $x$-dimension (streamwise) and represents an ensemble variable. A similar but slightly more general ensemble structure can be noted in the 2-D magnetohydrodynamic (Hartmann) flow \cite[eq. (59)]{xu2008stabilization}. But the 2-D PDEs in \cite{vazquez2007closed,xu2008stabilization} are not ensemble PDEs in the fully general sense because they do not exhibit coupling among the PDEs in the wavenumber $k$. 

The only 2-D PDE example tackled thus far with backstepping which exhibits a general ensemble structure is the spatially periodic flow on an annulus, called the ``convection loop,'' in which the thermal dynamics \cite[eq. (15)]{vazquez2006explicit}, which are boundary-controlled in the direction of the radial coordinate $r$, are coupled in the annular coordinate $\theta$ through the integral in $\phi$. However, unlike the hyperbolic ensembles in the present paper, for which we permit only a scalar input, an angle-dependent (i.e., function of space) boundary control is employed in \cite[eq. (31)]{vazquez2006explicit}. Finally, one can recognize a rudimentary form of a PDE ensemble system in the 1-D reaction–diffusions PDE with delayed distributed actuation \cite[eqs. (7)-(12)]{qi2019stabilization}, but this system also employs a functional input, parametrized by the ensemble variable $x$, rather than a scalar input like in the present paper.

The notion of ensemble systems has generated a rich literature on infinite-dimensional {\em non-PDE} systems. A few references that may stand out for being broadly known are \cite{li2009ensemble,li2010ensemble,helmke2014uniform,chen2019structure,chen2020ensemble,chen2021sparse,dirr2021uniform}. The inspiration for such systems comes from physics and other disciplines. This ODE-like ensemble control literature deals predominantly (though not exclusively) with controllability. It does not deal with stabilization. Given that ensemble systems are infinite-dimensional, with ensemble variables taking real (rather than integer) values, it is relevant to see their connection with PDE systems, as well as to recognize that, unlike PDEs, the operations that arise in the models of those ensemble systems exclude derivatives (while possibly including integrals) with respect to the ensemble variables. Hence, the research on control of ensemble PDEs is an important extension in two directions: an extension of PDE control to ensemble PDEs and an extension of the ODE-like ensemble control to PDE ensembles. 

\paragraph{Application inspiration}
The framework we introduce here, while primarily inspired by traffic and its continuum of classes of vehicles and drivers, is also applicable to fluids that are continuously stratified, structures that are not layered/laminated from discrete materials but manufactured with a continuous variation of materials, and many other applications that we expect to come to the attention of the research community as a result of this control design capability becoming known. 

While in this paper we focus on ensembles of coupled hyperbolic PDEs with boundary actuation, there is interest in this topic for other classes of PDEs, both for mathematical and application-driven reasons. One class is ensemble parabolic PDEs, which incorporates the dynamics of Lithium-ion batteries. Another class is coupled first-order hyperbolic PDEs with in-domain bilinear actuation, which incorporates numerous classes of population dynamics, ranging from biotechnology to epidemiology. This range of potential future developments underscores our interest in this topic and the relevance of the developments in this paper. 

\paragraph{Organization}
The class of ensemble hyperbolic PDEs is introduced in Section \ref{sec-model}. The PDE backstepping control design is presented in Section \ref{sec-design}. Closed-loop stability is shown in Section \ref{sec-stability}. Kernels that arise in the target system are derived in 
Section \ref{sec-cascade-kernel}. Kernels that arise in the backstepping transforms and the control law are studied in Section \ref{sec-kernel}. Simulation results are presented in Section \ref{sec-simulations}. Conclusions are provided in Section \ref{sec-conclusions}. 

\paragraph{Contribution}
The paper generalizes the result of \cite{dimeglio2013} from finitely many to infinitely (uncountably) many PDEs; from vector ensembles to functional ensembles. In addition to performing the necessary generalizations of all the technical results in \cite{dimeglio2013} from Euclidean to Hilbert spaces, there are unique technical challenges that we have to overcome in infinite dimension, such as the continuity of the characteristic curves with respect to the ensemble variable, which we tackle in Section \ref{subsec-cont-char-curves}.

\section{Coupled Hyperbolic PDE Ensemble with Boundary Actuation}
\label{sec-model}

We adapt the  results of Di Meglio's et al \cite{dimeglio2013} to the case where the homodirectional portion of the coupled  hyperbolic PDE system is no longer a finite-dimensional vector, $u(t,x) = [u_1(t,x), \ldots, u_n(t,x)]^{\rm T}$, with a spatial variable $x\in[0,1]$, but an uncountably infinite-dimensional vector, namely, an {\em ensemble}, denoted by $u(t,x,y)$, where $y\in[0,1]$. It is without loss of generality that both $x$ and $y$ are taken in the unit interval---any other finite interval can be shifted and scaled to the unit interval. (It should not be surprising that the unit interval is favored both for the spatial interval in the PDE literature and for the ensemble interval in the ensemble ODE control literature.)

We consider control design for an ensemble of coupled hyperbolic PDEs governed by
\begin{eqnarray}
   \label{osysuiota}
    u_t(t,x,y) + \lambda(x,y)u_x(t,x,y) &=& 
    \int_0^1 \theta(x,y,\eta) u(t,x,\eta)d\eta
    \nonumber \\ 
    &&+ W(x,y)v(t,x)
 \\
    \label{osysviota}
    v_t(t,x) - \mu(x)v_x(t,x) &=& \int_0^1\Xi(x,y)u(t,x,y)dy
    \nonumber\\ &&
\end{eqnarray}
with boundary conditions 
\begin{eqnarray}
\label{osysbound}
    u(t,0,y) &=& q(y)v(t,0), 
    \\
    v(t,1) &=& U(t)\,.
\end{eqnarray}

So, we proceed by considering a PDE with a state 
$u[t,x]  \in E$, where $E$ is the Hilbert space of square-integrable functions of $y\in[0,1]$ and $\langle p\,,\, q\rangle_E = \int_0^1 p(y) q(y)dy$ denotes $E$'s inner product.
In what follows, the notation $u[t,x]:= u(t,x,\cdot)$ serves to suppress the dependence on $y$ in $u(t,x,y)$, whereas $\{\cdot\}$ serves to denote an action of an operator on its argument (such an operator may be finite or infinite-dimensional, diagonal, integral, or differential).
We denote as $\L(E)$ the space of linear operators whose arguments are vectors in $E$ (i.e., functions of $y$). 

With this notation, we consider a PDE system on a domain with a spatial variable $x\in[0,1]$, with time $t\geq 0$, whose state $(u,v)$ has a component $u$ that takes values in $E$ for each $(t,x)$ and a component $v$ that takes real values for each $(t,x)$, and which is governed by the equations
\begin{eqnarray}
\label{osysu}
    u_t[t,x] + \boldsymbol{L}(x)\{u_x[t,x]\} &=& \boldsymbol{\Theta}(x)\{u[t,x]\}
        \nonumber \\ 
    &&+ W[x]v(t,x)
\\
\label{osysv}
    v_t(t,x) - \mu(x)v_x(t,x) &=& \langle \Xi[x]\; , \; u[t,x] \rangle_E
\end{eqnarray}
with boundary conditions 
\begin{eqnarray}
    \label{uboundq} 
    u[t,0] &=& q v(t,0)
    \\
    \label{ubound} 
    v(t,1) &=& U(t) 
\end{eqnarray}
where $q \in E$ and, for each $ x \in [0,1]$, $W[x],\,\Xi[x]$ are vectors in $ E$ (i.e., functions of $y$) and  $\boldsymbol{\Theta}(x),
\boldsymbol{L}(x) \in \L(E)$.
The integral operator $\boldsymbol{\Theta}(x)$ is defined, for $x \in [0,1]$, as
\begin{eqnarray}
\label{integral_op}
    \forall a \in E,\, \boldsymbol{\Theta}(x)\{a\}(y) = \int_0^1 \theta(x,y,\eta)a(\eta) d\eta
\end{eqnarray}
and, using $\theta(x,y,\eta) = \lambda(x,y)\delta(y-\eta)$, where $\delta$ is Dirac's delta function, the positive diagonal multiplication operator $\boldsymbol{L}$ is 
defined, as a special case of $\boldsymbol{\Theta}$, as
\begin{eqnarray}
    &\forall x \in [0,1],\,\forall a \in E, &\;  \boldsymbol{L}(x)\{a\}(y) = \lambda(x,y)a(y) \\ \nonumber
    &\text{with}\, &\forall y \in [0,1], \;  \lambda(x,y) > 0
\end{eqnarray}
For an integral operator defined as in \eqref{integral_op} we define its transpose operator $\boldsymbol{\Theta}^t(x)$ 
\begin{equation}
    \forall a \in E,\, \boldsymbol{\Theta}^t(x)\{a\}(y) = \int_0^1 \theta(x,\eta,y)a(\eta) d\eta
\end{equation}
We observe that it satisfies
\begin{equation}
   \forall a,\,b \in E^2, \langle a\,,\, \boldsymbol{\Theta}(x)\{b\} \rangle_E = \langle \boldsymbol{\Theta}^t(x)\{a\} \,,\, b \rangle_E
\end{equation}


\section{Backstepping Design}
\label{sec-design}

The target system is sought in the ``cascade'' form
\begin{eqnarray}
 \alpha_t[t,x] + \boldsymbol{L}(x)\{\alpha_x[t,x]\} &=& \boldsymbol{\Theta}(x)\{\alpha[t,x]\} + W[x]\beta(t,x) 
 \nonumber\\
 && + \int_{0}^{x}\kappa[x,\xi]\beta(t,\xi)d\xi \nonumber 
 \\ &&+ \int_{0}^{x}\boldsymbol{C}(x,\xi)\{\alpha[t,\xi]\}d\xi \label{alphasys}
 \\
    \beta_t(t,x) - \mu(x)\beta_x(t,x) &=& 0 \label{betasys}
\end{eqnarray}
with boundary conditions 
\begin{eqnarray}
    \label{boundtarget} 
    \label{bc-alpha}
    \alpha[t,0] &=& q v(t,0), 
    \\ \label{bc-beta}
    \beta(t,1) &=& 0\,,
\end{eqnarray}
where $q$ is the same arbitrary vector of $E$ as in \eqref{uboundq}, and where the kernel $\kappa[x,\xi] \in E$ and the operator $\boldsymbol{C}(x,\xi) \in \mathcal{L}(E)$ are yet to be specified, with $\T$ being the triangular domain
\begin{equation}
    \T = \{0\leq \xi\leq x\leq 1\}\,.
\end{equation}
With $k[x,\xi]$, which denotes a vector in $E$ (i.e., a function of $y$), and a scalar-valued $\Tilde{k}(x,\xi)$, both of which are functions of $x$ and $\xi$ in $\T$ yet to be specified, we introduce the backstepping transformation
\begin{eqnarray}
  \label{backalpha} \alpha[t,x]&=& u[t,x] \\
  \label{backbeta} \beta(t,x) &=& v(t,x) - \int_0^x \langle k[x,\xi]\,,\, u[t,\xi]\rangle_E d\xi     \nonumber \\ 
    &&- \int_0^x \Tilde{k}(x,\xi)v(t,\xi)d\xi
\end{eqnarray}
In order to deal with the joint systems $\gamma(t,x) = \left(\begin{array}{c}
\alpha[t,x]     \\    \beta(t,x)   \end{array}\right) $ and $w(t,x) = \left(\begin{array}{c}
u[t,x] \\v(t,x) \end{array} \right)$ we introduce $E^* = E \times \RR $ and $\langle \cdot, \, \cdot \rangle_*$ its canonical inner product.
We define 
\begin{equation}
\label{jointkernels-l}
k^*(x,\xi) =  \left(\begin{array}{c}
k[x,\xi]\\ \Tilde{k}(x,\xi)  \end{array} \right)
\end{equation}
Such that \eqref{backbeta} rewrites
\begin{eqnarray}
    \label{ktransf}
    && \beta(t,x) =  v(t,x) - \int_0^x \left \langle k^*(x,\xi)\,, \, w(t,\xi)  \right \rangle_*d\xi
\end{eqnarray}
We denote $\alpha[t] = \alpha(t,\cdot,\cdot)$ and introduce the $\L^2_{x,y}$ norm of the joint system: 
\begin{eqnarray}
    \left\| \left( \begin{array}{c}
      \alpha[t]\\
      \beta[t]
\end{array} \right) \right\|_{x,y} 
&=& \sqrt{\int_0^1 \|\alpha[t,x]\|_E^2 
+\beta(t,x)^2dx} 
\\ &=& \sqrt{\int_0^1 \|\gamma(t,x)\|_*^2dx}
\end{eqnarray}

Our goal is to find  conditions for $(k,\Tilde{k})$ to yield
\eqref{betasys}.
After some calculation (Appendix \ref{Appendix_kernel_design}) we derive a sufficient condition for this equality to be true: the kernels $k$ and $\Tilde{k}$  verifying for $(x,\xi) \in \mathcal{T} $ 
the ensemble (in $y$) coupled system of  first-order hyperbolic PDEs (in $x$ and $\xi$) given by
\begin{align}
\label{keq1}
    \mu(x)k_x[x,\xi] 
   - \boldsymbol{L}(\xi)\{k_\xi[x,\xi]\} 
   =&   \left( (\boldsymbol{L}'
   +\boldsymbol{\Theta}^t)(\xi)\right)\{k[x,\xi]\} \nonumber \\ &
   +\Xi[\xi] \Tilde{k}(x,\xi)
    \\
   \label{keq2}
    \mu(x) \Tilde k_x(x,\xi)
   +\mu(\xi) \Tilde k_\xi(x,\xi)
   =& -\mu'(\xi)\Tilde{k}(x,\xi) \nonumber \\ &
   +\langle W[\xi]\,,\, k[x,\xi]  \rangle_E  
\end{align}
with boundary conditions
\begin{align}
\label{kbound1}
(\boldsymbol{L}(x)+\mu(x)\boldsymbol{\rm{Id}})\{k[x,x]\} =& -\Xi[x]
\\ \label{kbound2}
\mu(0)\Tilde{k}(x,0) =& \langle q,  \boldsymbol{L}(0)\{k[x,0]\} \rangle_E
\end{align}
The existence of solutions $(k, \Tilde{k})$ of this system, which are regular enough, is for now assumed, but will be proved in Theorem \ref{well-posedness}.
For the reader's fuller clarity we point out that the ``scalar-valued'' form of the Goursat PDEs for $(k, \Tilde{k})$, in spatial variables $(x,\xi)$ and ensemble variable $y$, on the ``prysm-shaped'' (3D) domain $0\leq \xi\leq x\leq 1, \ 0\leq y \leq 1$, is given by
\begin{align}
\label{mkeq1}
    &\mu(x) \Tilde k_x(x,\xi)
   +\mu(\xi) \Tilde k_\xi(x,\xi)
   \nonumber\\
   & \quad = -\mu'(\xi)\Tilde k(x,\xi) 
   +\int_0^1 W(x,y) k(x,\xi,y) dy 
   \\
   \label{mkeq2}   &\mu(x)k_x(x,\xi,y)
   - \lambda(\xi,y) k_\xi(x,\xi,y)
   \nonumber\\
   &\quad =   
   \lambda_x(\xi,y)k(x,\xi,y)
   +
  \int_0^1 \theta(\xi,\eta,y,) k(x,\xi,\eta) d\eta
  \nonumber \\ &
   \quad\quad +\Xi(\xi,y) \Tilde k(x,\xi)
\end{align}
with boundary conditions
\begin{eqnarray}
\label{mkbound1}
k(x,x,y) &=& -\frac{\Xi(x,y)}{\lambda(x,y)+\mu(x)}
\\ \label{mkbound2}
\Tilde k(x,0) &=& \frac{1}{\mu(0)}\int_0^1 q(y)\lambda(0,y) k(x,0,y) dy,
\end{eqnarray}
where $k(x,\xi,\eta) = k[x,\xi](\eta)$.



We now formulate the main result about the controller design, similar to Theorem 3.2 in \cite{dimeglio2013}.

\begin{theorem}[Stability of the system]
\label{main_theorem}
Suppose that $W, \Theta \in \C^0([0,1],E)$, $\lambda,\mu \in \C^1([0,1],\RR^*_+)\times\C^1([0,1],\RR^*_+)$, and $\boldsymbol{\Theta} \in \C^0([0,1],\L(E))$ is uniformly bounded (in the sense of operator norm associated with $\|\cdot\|_E$). 
Then, 
the control law 
\begin{equation}
\label{feedback-law}
    U(t) = 
    \int_0^1 \left \langle k^*(1,\xi)\,, \, w(t,\xi)  \right \rangle_*d\xi
\end{equation}
where $k^*$
is a continuous solution of system \eqref{keq1}--\eqref{kbound2} as defined as in \eqref{ktransf}, guarantees that the closed-loop system  \eqref{osysu}--\eqref{ubound}, \eqref{feedback-law} is exponentially stable in the ensemble $\L^2_{x,y}$ norm given by $\left[\int_0^1\left(\int_0^1 u^2(t,x,y) dy + v^2(t,x)\right) dx\right]^{1/2}$.
\end{theorem}

It is helpful to note from \eqref{betasys}, \eqref{bc-beta} that the state $\beta$ of the autonomous system $\beta_t=\mu\beta_x, \beta(t,1)=0$ converges to zero in finite time, upon which, the remaining ensemble dynamics of the system $\alpha=u$ are governed, based on \eqref{backalpha}, \eqref{alphasys}, \eqref{bc-alpha}, by
\begin{eqnarray}
\label{eq-ensembtargetu}
u_t(t,x,y)  &=& -  \lambda(x,y)u_x(t,x,y)
\nonumber \\ 
&& +\int_0^1 \theta(x,y,\eta) u(t,x,\eta)d\eta 
\nonumber  \\ && 
+ \int_0^1\int_{0}^{x} c(x,\xi,y,\eta)u(t,\xi,\eta)\}d\xi d\eta
\quad \\
\label{eq-ensembtargetubound}
    u(t,0,y) &=&  0\,
\end{eqnarray}
where the kernel $c$ is defined such that 
\begin{align}
&\forall (x,\xi,a) \in \T \times E,\,\nonumber \\ 
&\qquad \boldsymbol{C}(x,\xi)\{a\}(y) = \int_0^1 c(x,\xi,y,\eta)a(\eta)d\eta
\end{align}
Note that every operator of $\L(E)$ does not always have such an integral expression in terms of a scalar kernel. However for this specific operator $ \boldsymbol{C}$ it is the case as we can derive this scalar expression with the kernel
\begin{equation}
    c(x,\xi,y,\eta) = W(x,y)k(x,\xi,\eta) + \int_\xi^x \kappa(x,s,y) k(s,\xi,\eta) ds\,.
\end{equation}
where $\kappa$ is the solution of the (spatial) integral equation
\begin{equation}
    \kappa(x,\xi,y) = W(x,y)\tilde k(x,\xi) + \int_\xi^x \tilde k(x,\xi)\kappa(x,s,y)  ds. 
\end{equation}
We return to these expression in the more compact form in \eqref{C_def}. 

The system \eqref{eq-ensembtargetu}, \eqref{eq-ensembtargetubound} is not necessarily finite-time stable but, thanks to the fact that the state $u$ convects only in the (increasing) direction of $x$, neither the integral in $\eta$ over $[0,1]$ nor the integral in $\xi$ over $[0,x]$ are destabilizing to the purely transport dynamics $u_t(x,y,t) + \lambda(x,y)u_x(x,y,t)=0$. 
In other words, this system has no feedback/recirculation from downstream (larger $x$) to upstream (smaller $x$), whereas the integration in $y$ represents only ``exchange'' of the state at a given transport position $x$, so, in addition to transport there may be growth but only over a finite distance. 
As, the system \eqref{eq-ensembtargetu}, \eqref{eq-ensembtargetubound} is exponentially stable, so is the system \eqref{alphasys}--\eqref{bc-alpha}, as (implicitly) asserted in Theorem \ref{main_theorem} and proved in the next section.

\section{Proof of Stability (Theorem~\ref{main_theorem})}
\label{sec-stability}

\begin{figure}[t]
    \centering
    \includegraphics[width=0.4\textwidth]{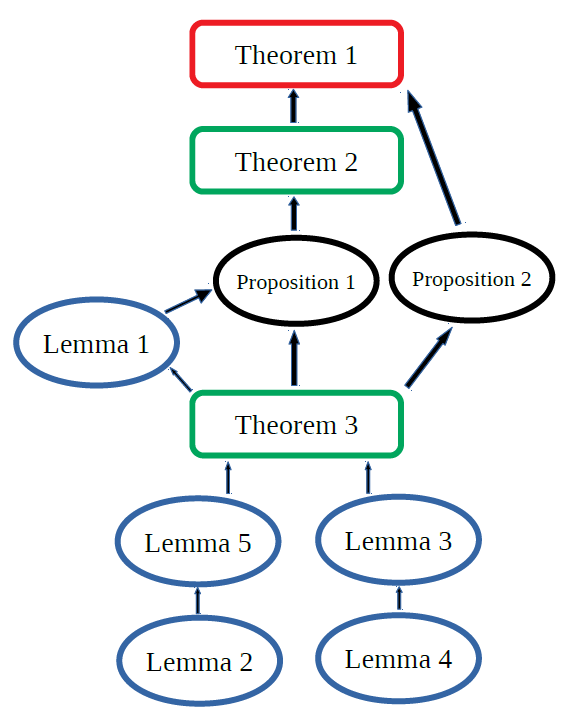}
    \caption{Interrelation among various technical and principal results in the paper, leading the to main result on feedback stabilization in Theorem \ref{main_theorem}. Theorem \ref{targetstab} establishes stability of the target system and Theorem \ref{well-posedness} establishes the well posedness of the kernels of the ensemble backstepping transformation.}
    \label{summary}
\end{figure}

The sketch of the proof is in two step: assess the stability of the target system and then show that the stability of both systems are equivalent. The relation between the different results developed is schematized in Figure \ref{summary}.

We first prove the exponential stability of the target system \eqref{alphasys}--\eqref{boundtarget}, similarly to what has been done in Lemma 3.1 of \cite{dimeglio2013}. We employ the following Lyapunov function, parametrized by $p>0$ and $\delta>0$, and given in the compact and scalar notations, respectively, by
\begin{eqnarray}
    V(t) &=& p \int_0^1 e^{-\delta x} \langle \alpha[t,x],\boldsymbol{L}(x)^{-1}\{\alpha[t,x]\}\rangle_E dx     \nonumber \\ 
    &&+ \int_0^1\frac{1+x}{\mu(x)}\beta(t,x)^2dx
    \\
    &=& p\int_0^1e^{-\delta x} \int_0^1 \frac{1}{\lambda(x,y)} \alpha^2(t,x,y) dy dx 
    \nonumber \\ 
    && + \int_0^1\frac{1+x}{\mu(x)}\beta(t,x)^2dx
    \label{Lyapounov}
\end{eqnarray}

\begin{theorem}[Target system stability]
Suppose that $x \mapsto  W[x]$ 
is a bounded functions (in the $\|\cdot\|_E$ norm), that $k^*$ is continuous 
and that both $\boldsymbol{L}^{-1}$ and $\boldsymbol{\Theta}$ are also bounded operator of $\L(E)$ for every $x$. Then the 0 equilibrium of the target system \eqref{alphasys}--\eqref{boundtarget} is exponentially stable in the $\L^2_{x,y}$ norm.
\label{targetstab}
\end{theorem}

One can check in Appendix \ref{target_system_stability} that \eqref{Lyapounov} is indeed an adequate Lyapunov function under the assumptions of Theorem \ref{targetstab}. 

We have introduced the continuous kernels $(k, \Tilde{k})$ (solution of system \eqref{keq1}--\eqref{keq2} with boundary conditions \eqref{kbound1}--\eqref{kbound2}) for transforming the original system into the target system. We also need to introduce the inverse kernel for transforming the target system back into the original system.

\begin{proposition}\label{invers_kernel}
There exist unique continuous kernels 
$l^* = \left(\begin{array}{c}
         l[x,\xi]\\ \Tilde{l}(x,\xi) \end{array} \right) \in \C([\T,E^*])$ 
such that the inverse of the backstepping transformation \eqref{backbeta} is given by 
\begin{equation}
v(t,x) = \beta(t,x) + \int_0^x \left \langle l^*(x,\xi)\,, \, \gamma(t,\xi)  \right \rangle_* d\xi
         \label{ltransf}
\end{equation}
\end{proposition}

The expressions of these kernels are derived in Appendix \ref{inverse_kernel_section}. Their existence, uniqueness, and continuity follow from the existence, uniqueness, and continuity of the direct kernels $k^*$ 

We now have all the necessary intermediate result to prove Theorem \ref{main_theorem}.

\vspace{1em}
\noindent
\\
\emph{Proof of Theorem \ref{main_theorem}} : Applying Theorem \ref{well-posedness} to the  system
\begin{eqnarray}
\label{Fgenericsystem1}\mu (x)F_x - \boldsymbol{L}(\xi)\{F_{\xi}\} &=& a[x,\xi]G+\boldsymbol{B}(x,\xi)\{F\}
\\
\label{Ggenericsystem1}\mu (x)G_x + \mu (\xi)G_{\xi} &=& d(x,\xi)G+ \langle e[x,\xi]\,,\,F \rangle_E
\\
F[x,x] &=& f[x]
\\
G(x,0) &=& \langle g[x]\,,\, F[x,0]\rangle_E \label{Ggenericsystem_bound1}
\end{eqnarray} 
with $a[x,\xi] =  \Xi[\xi],\, d(x,\xi) = -\mu'(\xi) ,\,e[x,\xi] = W[\xi] ,\,\boldsymbol{B}(x,\xi)= \boldsymbol{L'}(\xi) + \boldsymbol{\Theta}^t(\xi) ,\,f[x] = -(\boldsymbol{L}+\mu(x)\boldsymbol{Id})^{-1}\{\Xi[x]\},\, g[x] = \frac{1}{\mu(0)}\boldsymbol{L}(0)\{q\},\, G=\tilde{k},\, F=k$
we get the existence and continuity of the augmented kernel $k^*$. Following proposition \ref{invers_kernel} we also get the existence of its inverse kernel $l^*$.
Thanks to continuity of both kernels we introduce
\begin{eqnarray}
   K = \max_{(x,\xi) \in \T} \|k^*(x,\xi)\|_*, \, L = \max_{(x,\xi) \in \T} \|l^*(x,\xi)\|_*, \, 
\end{eqnarray}
Thanks to Theorem \ref{targetstab} we therefore have the existence of $C$ and $\varepsilon$ such that 
\begin{eqnarray}
   \|\gamma[t]\|_{x,y} \le C\|\gamma[0]\|_{x,y}e^{-\varepsilon t}
\end{eqnarray}
We then have, thanks to \eqref{ltransf}
\begin{eqnarray}
   \|w(t,x)\|_* &=& \left\|\left(\begin{array}{c}
\alpha[t,x]     \\    v(t,x)
\end{array}\right)\right\|_* 
=\bigg\|
\left(\begin{array}{c}
\alpha[t,x]     \\    \beta(t,x)\end{array}\right)   \nonumber \\
 &&+\left(\begin{array}{c}
0\\
 \int_0^x \left \langle l^*(x,\xi)\,, \, \gamma(t,\xi)  \right \rangle_* d\xi  \end{array}\right) \bigg\|_* 
\\
 &\le&  \|\gamma(t,x)\|_* \nonumber
 \\&&+ \left|\int_0^x \left \langle l^*(x,\xi)\,, \, \gamma(t,\xi)  \right \rangle_* d\xi \right|  
 \\ &\le& \|\gamma(t,x)\|_* + L\int_0^1 \left \| \gamma(t,\xi)  \right\|_* d\xi 
\end{eqnarray}
Therefore, using Jensen's inequality
 \begin{eqnarray}
    \|w[t]\|_{x,y}^2 &\le& \int_0^1 \|\gamma(t,x)\|_*^2dx \nonumber\\&+&L^2\left(\int_0^1 \left \| \gamma(t,\xi)  \right\|_* d\xi \right)^2 \nonumber
    \\&+&  \int_0^12L\|\gamma(t,x)\|_*\int_0^1 \left \| \gamma(t,\xi)  \right\|_* d\xi dx 
    \\&\le& (1+2L+L^2) \int_0^1 \|\gamma(t,x)\|_*^2dx
 \end{eqnarray}
Similarly, using \eqref{ktransf}
\begin{eqnarray}
   \|\gamma(t,x)\|_* &=& \left\|\left(\begin{array}{c}
u[t,x]     \\    \beta(t,x)
\end{array}\right)\right\|_* =\bigg\|
\left(\begin{array}{c}
u[t,x]     \\    v(t,x)\end{array}\right)  \nonumber
\\&&- \left(\begin{array}{c}
0\\
\int_0^x \left \langle k^*(x,\xi)\,, \, w(t,\xi)  \right \rangle_* d\xi  \end{array}\right) \bigg\|_*
\\
 &\le&  \|w(t,x)\|_* \nonumber\\&&+ \left|\int_0^x \left \langle k^*(x,\xi)\,, \, w(t,\xi)  \right \rangle_* d\xi \right|
 \\ &\le& \|w(t,x)\|_* + K\int_0^1 \left \| w(t,\xi)  \right\|_* d\xi 
\end{eqnarray}
which yields
\begin{equation}
 \|\gamma[t]\|_{x,y} \le (1+K)\|w[t]\|_{x,y}
\end{equation}
In the end we have 
\begin{eqnarray}
       \|w[t]\|_{x,y} &\le& (1+L)\|\gamma[t]\|_{x,y} \\&\le& (1+L)C\|\gamma[0]\|_{x,y}e^{-\varepsilon t}\\&\le& 
    (1+L)C(1+K)\|w[0]\|_{x,y}e^{-\varepsilon t}\hfill\blacksquare
\end{eqnarray}

\section{Existence and continuity of cascade connection Volterra kernel in the target system \eqref{alphasys}-\eqref{bc-beta}}
\label{sec-cascade-kernel}

We now turn our attention to the target system's equation \eqref{alphasys}-\eqref{betasys}, whose operator $\boldsymbol{C}$ and kernel $\kappa$ have not been specified yet. We assume the backstepping transformation kernels $(k,\Tilde{k})$ as known/given.
Plugging \eqref{backalpha} and \eqref{backbeta} into \eqref{osysv} gives
\begin{eqnarray}
0 &=& 
 \int_0^x \Bigg(\langle k[x,\xi]\,,\, u(t, \xi) \rangle_E W[x] - \boldsymbol{C}(x,\xi)\{u[t,\xi]\}  \nonumber
 \\&&+ \left. \int_0^{\xi} \langle k[x,s]\,,\,u(t,s)\rangle_E \kappa[x,\xi]ds \right) d\xi \nonumber
    \\ 
&&+ \int_0^x \Bigg(\Tilde{k}(x,\xi)W[x] - \kappa[x,\xi] \nonumber \\ &&+ \int_0^{\xi} \Tilde{k}(\xi,s)\kappa[x,\xi]ds \Bigg)v(t,\xi)d\xi  
\end{eqnarray}
This is true if, for all $((x,\xi),a) \in \T \times E$ , $\kappa$ and $\boldsymbol{C}$ satisfies
\begin{eqnarray}
\label{Volterra_kappa}
\kappa[x,\xi] &=& W[x]\Tilde{k}(x,\xi) + \int_{\xi}^x\Tilde{k}(s,\xi)\kappa[x,s]ds
\\
\boldsymbol{C}(x,\xi)\{a\} &=& \langle k[x,\xi]\,,\,a \rangle_E W[x] \nonumber \\ &&+   \int_{\xi}^x \langle k[s,\xi]\,,\, a\rangle_E\kappa[x,s]ds \label{C_def}
\end{eqnarray}
Equation \eqref{Volterra_kappa} is known as a Volterra equation of the second kind. It is known that such an equation admits an unique continuous solution in the scalar case (see eg \cite[ch. 8]{vrabieDifferentialEquationsIntroduction2016}). 
 $E$ being complete allows to use the method of successive approximations and to formulate the same results (see Lemma \ref{volterra_solution}) for infinite dimensional unknown $\kappa[x,\xi](y) = \kappa(x,\xi,y)$.
\begin{lemma}\label{volterra_solution}
If $\Tilde{k} \in \C(\T,\RR)$, $W \in \C([0,1],E)$, and $E$ is complete, the equation \eqref{Volterra_kappa} for the unknown $\kappa$ admits a continuous solution. 
Furthermore, this solution can be expressed in terms of a continuous inverse kernel $\tilde{k}^{\infty} \in \C(\T,\RR)$ as follows
\begin{align}
    \label{inverse_kernel}\kappa[x,\xi] =& \kappa_0[x,\xi] + \int_{\xi}^x\tilde{k}^{\infty}(s,\xi)\kappa_0[x,s]ds
    \\
    \text{with } \, \kappa_0[x,\xi] =& W[x]\Tilde{k}(x,s)
\end{align}
\end{lemma}

The existence and continuity of $\boldsymbol{C}$ readily follows from \eqref{C_def}, the continuity of $W$, which is assumed, and the continuity of $\kappa$ and $k$, which is proven in Lemma \ref{volterra_solution} and Theorem \ref{well-posedness}.

\begin{proposition}\label{C0_kappa_C}
The kernels $\kappa$ and $\boldsymbol{C}$ defined by \eqref{Volterra_kappa} and \eqref{C_def} exist and are continuous. 
\end{proposition}

\vspace{1em}
\noindent
\emph{Proof of Lemma \ref{volterra_solution} }:
Let $(\kappa_n)_n\in \mathcal{C}(\mathcal{T}, E)^{\NN}$ be a sequence of functions defined by induction \begin{eqnarray}
   \kappa_0: (x,\xi) \in \mathcal{T} &\mapsto& \Tilde{k}(x,s)W[x] \in E
   \\
   \kappa_{n+1}: (x,\xi) \in \mathcal{T} &\mapsto& 
   \kappa_0[x,\xi] \nonumber \\ &&+ \int_{\xi}^x\Tilde{k}(s,\xi)\kappa_n[x,s]ds 
\end{eqnarray}
We then define the functions $\Delta \kappa_{n+1} = \kappa_{n+1} - \kappa_n$, which verify the induction formula
\begin{eqnarray}
   \Delta \kappa_0[x,\xi] &=&\kappa_0[x,\xi]
  \\ \label{induc_kappa_def} \Delta \kappa_{n+1}[x,\xi] &=& \int_{\xi}^x\Tilde{k}[s,\xi]\Delta \kappa_{n}[x,s]ds
\end{eqnarray}
We introduce $M$ and $K$ such that \begin{eqnarray}
   \forall (x,\xi) \in \mathcal{T}, && |\Tilde{k}(x,\xi)| \le K
   \\ && \|W[x]\|_E \le M
\end{eqnarray}
We show by induction that \begin{equation}
\label{induction_kappa_c0_start}
    \forall (x,\xi) \in \mathcal{T},\;\|\Delta \kappa_n[x,\xi]\|_E \le MK \frac{\left(K(x-\xi)\right)^n}{n!}
\end{equation}
This is readily shown for $n= 0$. Suppose this inequality is true for some $n \in \NN$.
We then have for all $(x,\xi) \in \mathcal{T}$ \begin{eqnarray}
    \|\Delta \kappa_{n+1}[x,\xi]\|_E &\le& \int_{\xi}^x K\|\Delta \kappa_{n}[x,s] \|_Eds
    \\&\le&  KMK\frac{K^{n}}{n!}  \int_{\xi}^x (x-s)^n ds
    \\ &=& MK \frac{(K(x-\xi))^{n+1}}{(n+1)!} \label{induction_kappa_c0_end}
\end{eqnarray}
We can also show by induction that $\Delta \kappa_n$ are continuous functions of $x$ and $\xi$. The series $\sum \Delta \kappa_n$ being normal convergent in the complete space $\C (\T,E), \|\cdot\|_{\infty}$ we know that it uniformly converges to a continuous function $\kappa$. Furthermore $\kappa$ verifies \eqref{Volterra_kappa} thanks to telescoping. So far we have proven the existence of a continuous solution of \eqref{Volterra_kappa}. We would now like to show the existence of a continuous inverse kernel $\tilde{k}^{\infty}$ verifying \eqref{inverse_kernel}.

Similarly as what is done in \cite[2.3 p.19]{fedaabdelazizmustafaAnalyticalNumericalSolutions2014}, we define by induction the following sequence of kernels $(\Tilde{k}^n)_n\in \mathcal{C}(\mathcal{T}, \RR)^{\NN^*}$
\begin{eqnarray}
  \forall(x,\xi) \in \T,& \Tilde{k}^1(x,\xi) &= \Tilde{k}(x,\xi) 
  \\ &\Tilde{k}^{n+1}(x,\xi) &= \int_{\xi}^x \Tilde{k}(s,\xi) \Tilde{k}^n(x,s)ds\
\end{eqnarray}
Using \eqref{induc_kappa_def} and inverting the order of integration over the triangular domain of integration we can show by induction that
\begin{equation}
    \forall n \ge 1, \, \Delta \kappa_n[x,\xi] = \int_{\xi}^x \Tilde{k}^n(s,\xi)\kappa_0[x,s]ds
\end{equation}
We have, provided that the kernel series $\sum \Tilde{k}^n$ is uniformly convergent
\begin{eqnarray} 
   \kappa[x,\xi] &=& \sum_{n=0}^{\infty}\Delta\kappa_n[x,\xi] \\
   &=& \kappa_0[x,\xi] + \sum_{n=1}^{\infty}\Delta\kappa_n[x,\xi]
   \\
   &=& \kappa_0[x,\xi] + \sum_{n=1}^{\infty}\int_{\xi}^x \Tilde{k}^n(s,\xi)\kappa_0[x,s]ds
   \\ \label{intervertion}
   &=& \kappa_0[x,\xi] + \int_{\xi}^x\sum_{n=1}^{\infty} \Tilde{k}^n(s,\xi)\kappa_0[x,s]ds
\end{eqnarray}
We can show in a similar fashion than the induction \eqref{induction_kappa_c0_start}--\eqref{induction_kappa_c0_end} that
\begin{eqnarray}
\label{induction_k_c0}
   \forall n\ge 0, \,\forall (x,\xi) \in \T, \, |\Tilde{k}^{n+1}(x,\xi)| \le K \frac{(K(x-\xi))^n}{n!}
\end{eqnarray}
The series $\sum \Tilde{k}^n$ being normal convergent this allows us to undergo the intervertion made in \eqref{intervertion} and to introduce the continuous inverse kernel
\begin{eqnarray}
   \tilde{k}^{\infty}(s,\xi) = \sum_{n=1}^{\infty} \Tilde{k}^n(s,\xi) \label{inverse_kernel_def}
\end{eqnarray}
Combining \eqref{intervertion} and \eqref{inverse_kernel_def} we get \eqref{inverse_kernel} which concludes the proof. \hfill$\blacksquare$

\section{Existence and Continuity of  of solutions to Kernel PDE's}
\label{sec-kernel}

In order to prove the existence of the different Kernels used (such as $k,\Tilde{k}$) we adapt the Theorem 5.3 of \cite{dimeglio2013} into theorem \ref{well-posedness}. We are interested in the generic system defined on $\mathcal{T} = \{(x,\xi), 0\le \xi \le x \le 1\}$
\begin{eqnarray}
\label{Fgenericsystem}\mu (x)F_x - \boldsymbol{L}(\xi)\{F_{\xi}\} &=& a[x,\xi]G+\boldsymbol{B}(x,\xi)\{F\}
\\
\label{Ggenericsystem}\mu (x)G_x + \mu (\xi)G_{\xi} &=& d(x,\xi)G+ \langle e[x,\xi]\,,\,F \rangle_E
\end{eqnarray}
with boundary conditions
\begin{eqnarray}
F[x,x] &=& f[x]
\\
G(x,0) &=& \langle g[x]\,,\, F[x,0]\rangle_E \label{Ggenericsystem_bound}
\end{eqnarray}
One should note that $G(x,\xi)$, $\mu(x)$ and $d(x,\xi)$ are scalars, $F[x,\xi]$, $a[x,\xi]$, $e[x,\xi]$, $f[x]$ and $g[x]$ are vectors of $E$ while $\boldsymbol{B}(x,\xi)$ and $\boldsymbol{L}(x,\xi)$ are operators of $\mathcal{L}(E)$. Furthermore we also want $\mu$ to be a strictly positive valued function and $\boldsymbol{L}(x,\xi)$ to be a positive multiplication operator associated to $\lambda$: $\forall a \in E, \;  \boldsymbol{L}(x)\{a\}(y) = \lambda(y, x)a(y)$

\begin{theorem}[Kernel well posedeness] 
\label{well-posedness}
Assume that $a,e \in \C^0(\T, E)$, $d \in \C^0(\T, \RR)$, $f,g \in \C^0([0,1], E)$, $\boldsymbol{B} \in \C^0(\T, \L(E))$, $\mu \in \C^1([0,1], \RR^*_+)$, $\lambda \in \C^1([0,1]^2, \RR^*_+)$. Then there exists a continuous solution to system \eqref{Fgenericsystem}--\eqref{Ggenericsystem_bound}.
\end{theorem}

The sketch of the proof is as follows: transform PDE equations \eqref{Fgenericsystem}-\eqref{Ggenericsystem} into ODE integral equations thanks to integrating the PDEs along well-designed characteristic curves; then find a solution of these equations using the method of successive approximations.
\subsection{Characteristic curves}
\label{sec-charac_curves}
Let us fix $(x,\xi) \in \mathcal{T}$ for the following of this section, in which we will introduce their respective characteristic curves.
To define those curves we first need to present the two following lemmas : Lemmas 5.1 and Lemma 5.2 from \cite{dimeglio2013} .
\begin{lemma}
\label{5.1}
Let $(w_0,z_0)\in \RR^2$ be such that $0 \le w_0 \le z_0\le1$  
and $ h \in \mathcal{C}^1([0,1])$ be such that $\forall x \in [0,1], \, h(x)<0$ . 
Then, if $w$ and $z$ are the maximal solutions of the following Cauchy
problems
\begin{equation}
    w'(s) = h(w(s)),\, w(0) = w_0, \, z'(s) = h(z(s)), \, z(0) = z_0
\end{equation} 
then, there exists $T>0$ such that $w(T) = 0$ and $z(T) \ge 0$.
\end{lemma}
\begin{lemma}
\label{5.2}
Let $(w_0,z_0)\in \RR^2$ be such that $0 \le w_0 \le z_0\le1$ 
and $ h,g \in \mathcal{C}^1([0,1])$ be such that $\forall x \in [0,1],\, h(x)<0,\, g(x) >0$ . 
Then, if $w$ and $z$ are the maximal solutions of the following Cauchy
problems
\begin{equation}
w'(s) = g(w(s)), \, w(0) = w_0, \, z'(s) = h(z(s)), \,z(0) = z_0
\end{equation} 
then there exists $T>0$ such that $w(T) = z(T)$.
\end{lemma}
Applying Lemma \ref{5.1} to 
\begin{equation}
h = -\mu, \, w_0 = \xi, \, z_0 = x 
\end{equation} and defining 
\begin{eqnarray}
&&s_F = T \label{s_Fdef}\\ &&s \in [0, s_F] \mapsto \chi(s) = z(T-s), \,\chi_0 = z(T)
\\ &&s \in [0, s_F] \mapsto \zeta(s) = w(T-s)
\end{eqnarray}
we have the existence of the curve $(\chi, \zeta)$ and $s_F$ such that
\begin{eqnarray}
  \dot{\chi}(s) &=& \mu(\chi(s)) \\ \chi(s_F) &=& x, \,\chi(0) = \chi_0 \ge 0
 \\
\nonumber\\
  \dot{\zeta}(s) &=& \mu(\zeta(s))  \\ \zeta(s_F) &=& \xi, \,\zeta(0) =  0
\end{eqnarray}
Similarly, for each $y$, we apply Lemma \ref{5.2} to $w_{y}$ and $z$, respective solutions of the following Cauchy problems
\begin{eqnarray}
\dot{z} &=& h(z),  \, z_0 = x 
\\ \dot{y} &=& g_{y}(y) , \, w_0 = \xi
\\ \text{where} \qquad h &=& -\mu,\, g_{y} = \lambda(y,\cdot)
\end{eqnarray}
If we define 
\begin{eqnarray}
&&s_f(y) = T_{y} \label{s_f_def}
\\&&\hat{x}_0(y) = w_{y}(T_{y}) = z(T_{y})
\\ &&s \in [0, s_F] \mapsto \hat{\xi}(y,s) = w_{y}(T_{y}-s), 
\\ &&s \in [0, s_f(y)] \mapsto \hat{x}(y,s) = z(T_{y}-s)
\end{eqnarray}
we have, for each $y$, a characteristic curve $(\hat{x}(y, \cdot),\hat{\xi}(y, \cdot))$ and $s_f(y)$ such that
\begin{eqnarray}
  \frac{d\hat{x}}{ds}(y,s) &=& \mu(\hat{x}(y,s))  \label{xiotaeq}\\
  \hat{x}(y,s_f(y)) &=& x, \,\hat{x}(y,0) = \hat{x}_0(y) \label{xiotabound}
 \\
\nonumber\\
 \frac{d\hat{\xi}}{ds}(y,s) &=& -\lambda(y,\hat{\xi}(y,s))  \label{xiiotaeq} \\  \hat{\xi}(y,s_f(y)) &=& \xi, \,\hat{\xi}(y,0) = \hat{x}(y,0) \label{xiiotabound}
 \end{eqnarray}
Integrating \eqref{Fgenericsystem} for a given $y$ (we are forced to reason coordinates-wise because Lemma \ref{5.2} is only a scalar result) over the $(\hat{x}(y,\cdot),\hat{\xi}(y,\cdot))$ curve gives us
\begin{align}
 F(x,\xi,y) =& f(\hat{x}_0(y),y)  \nonumber
 \\&+ \int_0^{s_f(y)} \Big[a(\hat{x}(y,s), \hat{\xi}(y, s), y)G(\hat{x}(y,s), \hat{\xi}(y, s))
 \nonumber \\& +\boldsymbol{B}(\hat{x}(y,s), \hat{\xi}(y, s))\{F[\hat{x}(y,s), \hat{\xi}(y, s)]\}(y) \Big]ds
\end{align}
Similarly  with \eqref{Ggenericsystem} we have
\begin{align}
   \label{GasODE}
 G(x,\xi) =& \langle g[\chi_0], F[\chi_0,0] \rangle_E  \nonumber \\ &+ \int_0^{s_F} \Big[d(\chi(s), \zeta(s))G(\chi(s),\zeta(s))
 \nonumber \\& + \Big\langle e[\chi(s), \zeta(s)]\,, \, F[\chi(s), \zeta(s)] \Big\rangle_E \Big] ds
\end{align}

\subsection{Continuity of characteristic curves}
\label{subsec-cont-char-curves}

{Thus far in the paper most of the material presented can be regarded as generalization of the respective material in \cite{dimeglio2013} from Euclidean to Hilbert spaces. In this section we  tackle a problem that is unique to the \emph{ensemble} PDE and has no counterpart in the finite dimensional case, namely the {\em continuity} of the characteristic curves $\hat{x},\hat{\xi}$ and their endpoint's abscissa $s_f, s_F$ with respect to $y$.}
If we want to use the method of successive approximations we need for $F_n[x,\xi]$ to be in $E$, i.e., to have a certain regularity with respect to $y$. We therefore need $s_f(y)$ to have a certain regularity with respect to $y$, which will be proven in Lemma \ref{s_Fc0}. This problem is specific to the infinite dimensional case as a finite collection of scalars can always be seen as a vector of the corresponding Euclidean space whereas a continuously infinite collection is not necessarily a vector of E.

\begin{figure}[t]
    \centering
    \includegraphics[width=0.48\textwidth]{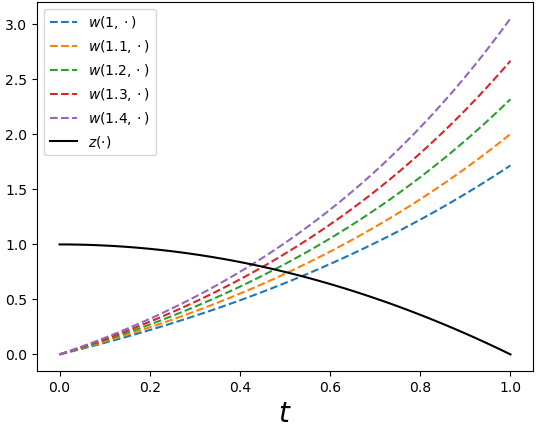} 
    \caption{Visualization of the continuity of $s_f(y,x,\xi)$ regarding the parameter $y$ with fixed $(x,\xi)$. This point is the abscissa at which the decreasing function $t \mapsto z(t)$ intersects the increasing function $t \mapsto w(y,t)$,  continuously parameterized by $y$.}
    \label{fig:sf_c0}
\end{figure}

Contrary to the previous subsection \ref{sec-charac_curves} we no longer treat $x$ and $\xi$ as fixed. Every object defined in the previous subsection now potentially varies with ($x,\xi)$. In particular, $s_f(y)$ becomes $s_f(y,x,\xi)$ and $s_F$ becomes $s_F(x,\xi)$.
Note that we do not explicitly write the dependence of the characteristic curves $z,w_y$ on their initial conditions $x,\xi$. We do so for clarity. One should only remember that this dependence of the characteristic curves on their initial conditions is continuous provided that $\lambda(x,y)$ and $\mu(x)$ are Lipschitz. The Lipschitzness property is also enough to have $(y,x,\xi) \mapsto s_f(y,x,\xi)$ continuous as stated in Lemma \ref{sf_co}. The statement of this lemma is illustrated in Figure \ref{fig:sf_c0}.
Regularity of $\hat{\xi}(y,\cdot)$ and $\hat{x}(y,\cdot)$ with respect to $(y, x, \xi)$ is then ensured as a composition of continuous functions.
Similarly for $\zeta$ and $\chi$; those curves depend continuously on $(x,\xi)$ if  $(x, \xi) \in  \T \mapsto s_F(x,\xi)$ is continuous (see Lemma \ref{s_Fc0}) and $x \mapsto \mu(x)$ is Lipschitz.
\begin{lemma}
\label{sf_co}
Assume that $(x, y) \mapsto \lambda(x,y)$ and $x\mapsto \mu(x)$ are respectively $L_{\lambda}$-Lipschitz and $L_{\mu}$-Lipschitz.
We then have that $(y, x, \xi) \in [0,1] \times \T \mapsto s_f(y,x,\xi)$, the function of $y$ implicitly defined in \eqref{s_f_def}, is Lipschitz.
\end{lemma}

\emph{Proof}:
By integrating \eqref{xiotaeq}--\eqref{xiiotaeq} and using both boundary conditions \eqref{xiiotabound}, \eqref{xiotabound} we have, for any $(x,\xi) \in\T$
\begin{align}
    x - \xi =& \int_0^{s_f(y,x,\xi)} \mu(\hat{x}(y,s)) + \lambda(y,\hat{\xi}(y,s)) ds
    \\
\label{sf_def}
    x - \xi =& \int_0^{s_f(y,x,\xi)} -h(z(s)) + g_{y}(w_{y}(s)) ds
\end{align}
We introduce $\varepsilon_1$ such that $ \forall (y, x, \xi) \in [0,1]^3,$
\begin{equation}
\label{epsilon1}
    \mu(x) + \lambda(y, \xi) > \varepsilon_1
\end{equation}
Using \eqref{epsilon1} in \eqref{sf_def} gives
\begin{equation}
\label{M_s}
      0 \le \xi- x + s_f(y,x,\xi) \varepsilon_1
\end{equation}
which implies that \begin{equation}
    M_s = \frac{1}{\varepsilon_1} \ge \frac{x-\xi}{\varepsilon_1}
\end{equation} is a uniform bound of the function $(y, x, \xi) \mapsto s_f(y, x, \xi) $.
\\
We consider the augmented state $q(t) = \left(\begin{array}{c}
        y\\
       z(t)
       \\ w_{y}(t)
   \end{array} \right)$. Its dynamic is
\begin{equation}
\label{augmented_dynamic}
   \dot{q}(t) = \left(\begin{array}{c} 0\\-\mu(z)\\ \lambda(y,w_{y})
   \end{array} \right)
\end{equation}
Let $q_1^0 = \left(\begin{array}{c}
      y_1 \\ x_1 \\ \xi_1 
\end{array}\right) \in [0,1] \times \T$ and $q_2^0 = \left(\begin{array}{c}
      y_2 \\ x_2 \\ \xi_2 
\end{array}\right) \in [0,1] \times \T$.
We introduce $q_1(t) = \left(\begin{array}{c}
      y_1 \\ z_1(t) \\ w_1(t)
\end{array}\right)$ and $q_2(t)  = \left(\begin{array}{c}
      y_2 \\ z_2(t) \\ w_2(t)
\end{array}\right)$, respectives solutions of the system \eqref{augmented_dynamic} with initial conditions $q_1^0, q_2^0$.
Thanks to Lipschitz continuity of $\mu$ and $\lambda$ we know that there exists $M$ and $D$ such that 
\begin{eqnarray}
\label{lipschitz_augmented} 
   &\forall t>0,\;\|q_1(t) - q_2(t)\| \le Me^{Dt} \|q_1^0 - q_2^0\|
\end{eqnarray}
Using \eqref{sf_def} we get 
\begin{eqnarray}
 0 &=& \xi_1-\xi_2 -(x_1-x_2) \nonumber \\&&+ \int_0^{s_f(y_1,x_1,\xi_1)} -h(z_1(s)) +
 g_{y_1}(w_{1}(s))ds
\nonumber \\
&& - \int_0^{s_f(y_2,x_2,\xi_2)} -h(z_2(s)) + g_{y_2}(w_{2}(s))ds
\\
&=& \xi_1-\xi_2 -(x_1-x_2)\nonumber \\&&+ \int_{s_f(y_2,x_2,\xi_2)}^{s_f(y_1,x_1,\xi_1)} -h(z_1(s)) + g_{y_1}(w_{1}(s))ds \nonumber
\\
&&+ \int_0^{s_f(y_2,x_2,\xi_2)} \Big(-h(z_1(s)) + g_{y_1}(w_{1}(s))  \nonumber \\ &&
+h(z_2(s)) - g_{y_2}(w_{2}(s))\Big)ds \label{continuous_point}
\end{eqnarray}
Finally, by integrating \eqref{epsilon1} between $s_f(y_1,x_1,\xi_1))$ and $s_f(y_2,x_2,\xi_2)$ and using  \eqref{continuous_point} we get
\begin{eqnarray}
  &&\varepsilon_1|s_f(y_1,x_1,\xi_1)) - s_f(y_2,x_2,\xi_2)| \\
  &\le& \left|\int_{s_f(y_2,x_2,\xi_2)}^{s_f(y_1,x_1,\xi_1)} g_{y_1}(w_{1}(s))-h(z_1(s))ds   \right|
  \\&=& 
  \Bigg| \xi_1-\xi_2 -(x_1-x_2) 
  \\ &&+\int_0^{s_f(y_2,x_2,\xi_2)} g_{y_1}(w_{1}(s)) - g_{y_2}(w_{2}(s))ds \nonumber \\&& +\int_0^{s_f(y_2,x_2,\xi_2)} h(z_2(s))-h(z_1(s))ds \Bigg|
  \\ &\le& \sqrt{2} \|q_1^0 - q_2^0\| + L_{\lambda}\int_0^{M_s} \left\| \left( \begin{array}{c}
       y_1 - y_2\\ w_1(s) - w_2(s)
  \end{array}\right) \right\| ds\nonumber \\&& + L_{\mu}\int_0^{M_s}|z_1(s) - z_2(s)| ds
  \\ &\le&  
  \sqrt{2} \left( M_s(L_{\mu}+L_{\lambda})Me^{DM_s}+1\right) \|q_1^0 - q_2^0\| 
\end{eqnarray}
The last inequality is obtained using \eqref{lipschitz_augmented} and concludes the proof.\hfill $\blacksquare$

\begin{corollary}
\label{cor,x_0}
Under the same assumption as in Lemma \ref{sf_co} the implicit function  $(y, x,\xi) \in [0,1] \times \T \mapsto \hat{x}_0(y,x,\xi)$ is continuous.
\end{corollary}

\emph{Proof}: Integrating \eqref{xiotaeq} between $0$ and $s_f(y, x,\xi)$ gives
\begin{equation}
    \hat{x}_0(y,x,\xi) = x + \int_0^{s_f(y,x,\xi)}h(z(s))ds
\end{equation}
From this expression one can easily deduce the continuity of $\hat{x}_0$.
We could also use $\hat{x}_0(y,x,\xi) = z(s_f(y,x,\xi))$. \hfill $\blacksquare$
\\
\begin{lemma}
\label{s_Fc0}
Assume that $x\mapsto \mu(x)$ is positive valued and $\C^0$.
We then have that $(x, \xi) \in  \T \mapsto s_F(x,\xi)$, the implicit function of $(x, \xi)$ defined in \eqref{s_Fdef}, is $\C^1$.
\end{lemma}

\emph{Proof}: 
We define 
\begin{equation}
    \forall s,\xi \in \RR\times [0,1],\qquad\Phi(s,\xi) = \int_0^s \mu(\chi(\eta)) d\eta - \xi
\end{equation}
Let $\varepsilon$ be a strictly positive lower bound of $\mu$. Since $\Phi(0,\xi) \le 0$ and $\Phi(\frac{1}{\varepsilon},\xi) \ge 0$, applying the intermediate value Theorem to the continuous function $\Phi(\cdot,\xi)$ we get that for each $\xi \in [0,1]$ there exists $s_F(\xi) \in [0,\frac{1}{\varepsilon}]$ such that $\Phi(s_F(\xi),\xi) = 0$. 
Note that $\Phi$ is $\C^1$ and that $ \frac{\partial \Phi}{\partial s} = \mu(\chi(s)) > 0 $, therefore using the implicit function Theorem we have that $\xi \mapsto s_F(\xi)$ is locally $\C^1$. Since this function is globally defined we know that it is globally $\C^1$.  $(x, \xi) \in  \T \mapsto s_F(x,\xi) = s_F(\xi)$ is then readily continuous. \hfill $\blacksquare$

Following this result we can show a similar result as Corollary \ref{cor,x_0}.

\begin{corollary}
\label{cor,Xi_0}
Under the  assumption in Lemma \ref{s_Fc0}, the implicit function $(x,\xi) \mapsto \chi_0(x,\xi)$ is continuous.
\end{corollary}

\subsection{Successive approximation}
We know have all the tools to prove Theorem \ref{well-posedness}.

\emph{Proof}: We define by induction the following sequence of functions ($G$ is scalar valued and $F$ is a vector of $E$), initialized by the null function
\begin{align}
&F_{n+1}(y,x,\xi) - f(\hat{x}_0(x,\xi,y),y) \nonumber \\ =&\int_0^{s_f(y,x,\xi)} \Big[a(\hat{x}(y,s), \hat{\xi}(y, s), y)G_n(\hat{x}(y,s), \hat{\xi}(y, s))
\nonumber \\& +\boldsymbol{B}(\hat{x}(y,s), \hat{\xi}(y, s))\{F_n[\hat{x}(y,s)\,,\, \hat{\xi}(y, s)]\}(y)\Big]ds \label{equ-F-induction}
\end{align}
\begin{align}
G_{n+1}(x,\xi) =& \langle g[\chi_0(x,\xi)]\,,\, F_{n+1}[\chi_0(x,\xi),0] \rangle_E  \nonumber \\&+ \int_0^{s_F(x,\xi)} \Big[d(\chi(s), \zeta(s))G_n(\chi(s), \zeta(s))
\nonumber \\& + \Big\langle e[\chi(s), \zeta(s)]\,,\, F_n[\chi(s),\zeta(s)] \Big\rangle_E\Big] ds
\end{align}
 For clarity and rigor purposes we use the scalar notation to formulate those results but one could see the infinitely many equations \eqref{equ-F-induction} as one integration along the curve of $E$ $s \mapsto 
 (\hat{x}[s], \hat{\xi}[s])$.
 We can show by induction that these functions are continuous regarding $(x,\xi)$ as they are integrals over continuous bounds of continuous functions. The continuity of the bounds   $s_f(y,x,\xi)$ and $s_F(x,\xi)$ are ensured by Lemma \ref{sf_co} and Lemma \ref{s_Fc0}. The continuity of initial conditions $\hat{x}_0, \chi_0$ is ensured by Corollaries \ref{cor,x_0} and \ref{cor,Xi_0}. Similarly, every $F_n$ are continuous regarding $y$ which ensures that $F_n[x,\xi]$ are indeed vectors of $E$.
 We introduce the different quantities $\|\boldsymbol{B}\|, \|a\|, \|e\|,|d|, \|f\|, M_{\mu}, \|g\|, M_{\lambda}$ such that 
 \begin{align}
    \forall (x,\xi) \in \mathcal{T}, 
    \;&\forall u \in E, \|\boldsymbol{B}(x,\xi)\{u\}\|_E\le \|\boldsymbol{B}\| \|u\|_E
    \\ & \|a(x,\xi)\|_E \le \|a\|,\; \|e(x,\xi)\|_E \le \|e\|
    \\& |d(x,\xi)| \le |d|
    \\& \frac{1}{|\mu(x)|} \le M_{\mu}
    \\& \|f(x)\|_E \le \|f\|, \; \|g(x)\|_E \le \|g\|
    \\& \forall y \in [0,1], \; \frac{1}{|\lambda(y,x)|} \le M_{\lambda}
 \end{align}
 We define
 \begin{align}
    \forall n \ge 0, &\forall (x,\xi) \in \T, \nonumber \\& \Delta F_n(x,\xi,y) =  F_{n+1}(x,\xi,y) - F_n(x,\xi,y)
    \\
    &\Delta G_n(x,\xi) = G_{n+1}(x,\xi) - G_n(x,\xi)
 \end{align}
By drawing inspiration from the proof of \ref{volterra_solution}, we can again find $M$ and $K$ big enough (namely $K \ge M_\lambda (\|g\|+1)(\|a\| +\|\boldsymbol{B}\|)
  +
 M_\mu (\|e\|+|d|)
 $ and $M \ge (\|g\|+1)\|f\|$) such that we have, by induction 
 \begin{eqnarray}
    \forall n \ge 0, \forall x,\xi \in \T, \|\Delta F_n[x,\xi]\|_E &\le& M\frac{(Kx)^n}{n!}\label{F_induction}
    \\
    |\Delta G_n(x,\xi)| &\le& M\frac{(Kx)^n}{n!} \label{G_induction}
 \end{eqnarray}
In order to prove this induction result we first states the following inequalities
\begin{eqnarray}
\int_0^{s_f(y,x,\xi)} \hat{x}(y,s)^m ds \le M_\lambda \frac{x^{m+1}}{m+1} \label{ineqxhat}
\\
\int_0^{s_F(x,\xi)} \chi(s)^m ds \le M_\mu \frac{x^{m+1}}{m+1} \label{ineqchi}
\end{eqnarray}
Indeed using the change of variable $\varsigma = \hat{x}(y,s)$ and noticing $ d\varsigma = \lambda(y,x)ds $ we get
\begin{eqnarray}
&&\int_0^{s_f(y,x,\xi)} \hat{x}(y,s)^m ds = \int_{\hat{x}_0(y,x,\xi)}^{x} \frac{\varsigma^m}{\lambda(y,x)} d\varsigma 
\nonumber\\ && 
\le \int_{0}^{x} M_\lambda \varsigma^m d\varsigma = M_\lambda \frac{x^{m+1}}{m+1}
\end{eqnarray}
\eqref{ineqchi} is obtained in a similar fashion with the change of variable $\varsigma = \chi(s)$.
Lastly note that 
\begin{equation}
    \chi_0(x,\xi) = x - \int_0^{s_F(x,\xi)} \mu(\chi(s))ds \le x \label{Chi0<x}
\end{equation}
The inequalities \eqref{F_induction}, \eqref{G_induction} hold for $n =0$
\begin{align}
&\|\Delta F_0[x,\xi]\|_E = \|F_1[x,\xi]\|_E =  \|f[\hat{x}_0(y,x,\xi)]\|_E \le M
\\
&|\Delta G_0(x,\xi)| \le \|g[\chi_0]\|_E \|f(\chi_0)\|_E \le M
\end{align}
For $n \in \NN$, assume that \eqref{F_induction}, \eqref{G_induction} are true.
For all $(x,\xi)$ in $\T$, using \eqref{ineqxhat} we have 
\begin{align}
&\|\Delta F_{n+1}[x,\xi]\|_E \nonumber \\ =& \bigg\|\int_0^{s_f(y,x,\xi)} \Big[a[\hat{x}(y,s), \hat{\xi}(y, s)]\Delta G_n(\hat{x}(y,s), \hat{\xi}(y, s))
 \nonumber \\& +\boldsymbol{B}(\hat{x}(y,s), \hat{\xi}(y, s))\{\Delta F_n[\hat{x}(y,s)\,,\, \hat{\xi}(y, s)]\}\Big]ds  \bigg\|_E
 \\
 \le& \int_0^{s_f(y,x,\xi)} \Big\|a[\hat{x}(y,s), \hat{\xi}(y, s)]\Delta G_n(\hat{x}(y,s), \hat{\xi}(y, s)) \Big\|_Eds 
 \nonumber \\+&\int_0^{s_f(y,x,\xi)}  \Big\|\boldsymbol{B}(\hat{x}(y,s), \hat{\xi}(y, s))\{\Delta F_n[\hat{x}(y,s)\,,\, \hat{\xi}(y, s)]\}\Big\|_Eds  
 \\
 \le& \|a\|\int_0^{s_f(y,x,\xi)}M\frac{(K\hat{x}(y,s))^n}{n!}ds \nonumber\\ &+\|\boldsymbol{B}\|\int_0^{s_f(y,x,\xi)}M\frac{(K\hat{x}(y,s))^n}{n!}ds
 \\ \le& M\frac{K^n}{n!}(\|a\| +\|\boldsymbol{B}\|)\int_0^{s_f(y,x,\xi)} \hat{x}(y,s)^m ds \\ \le& M\frac{K^n}{n!}(\|a\| +\|\boldsymbol{B}\|)
 M_\lambda \frac{x^{n+1}}{n+1} \label{Knxn+1}
 \\ \le& M\frac{(Kx)^{n+1}}{(n+1)!}
\end{align}
and, using \eqref{ineqchi}-\eqref{Chi0<x}, \eqref{Knxn+1}
\begin{align}
  |\Delta G_{n+1}(x,\xi)| =& \bigg|\langle g[\chi_0(x,\xi)]\,,\, \Delta F_{n+1}[\chi_0(x,\xi),0] \rangle_E  \nonumber \\ &+ \int_0^{s_F(x,\xi)} \Big[d(\chi(s), \zeta(s))\Delta G_n(\chi(s), \zeta(s))
 \nonumber \\& + \Big\langle e[\chi(s), \zeta(s)]\,,\, \Delta F_{n}[\chi(s),\zeta(s)] \Big\rangle_E\Big] ds\bigg|
 \\ \le& \|g\|\|\Delta F_{n+1}[\chi_0(x,\xi),0]\|_E  \nonumber \\&+ \int_0^{s_F(x,\xi)} |d| |\Delta G_n(\chi(s), \zeta(s))|ds
 \nonumber \\& + \int_0^{s_F(x,\xi)} \|e\| \|\Delta F_{n}[\chi(s),\zeta(s)]\|_E ds
 \\ \le& 
 \|g\|M\frac{K^n\chi_0^{n+1}}{(n+1)!}(\|a\| +\|\boldsymbol{B}\|)
 M_\lambda  \nonumber \\&+ (\|e\|+|d|)\int_0^{s_F(x,\xi)}M\frac{(K\chi(s))^n}{n!}ds
 \\ \le& M\frac{K^n x^{n+1}}{(n+1)!}(M_\lambda \|g\|(\|a\| +\|\boldsymbol{B}\|)
  \nonumber \\& +
 M_\mu (\|e\|+|d|)
 )
 \\ \le& M\frac{(Kx)^{n+1}}{(n+1)!}
 \end{align}
This ensures the normal convergence of the sequences (viewed as the series of their respective differences $\Delta$) and henceforth their continuity similarly as what 
done in Lemma \ref{volterra_solution}.\hfill$\blacksquare$

\begin{figure}[p]
\centering
\begin{subfigure}[b]{\linewidth}
    \centering
    \includegraphics[width=.9\textwidth]{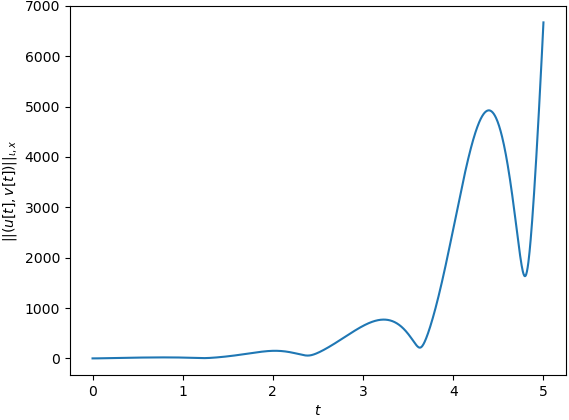} 
    \caption{Open loop. Instability results from the signal amplification due to the coupling of $u$ and $v$.}
    \label{open_loop_subfig}
\end{subfigure}
\begin{subfigure}[b]{\linewidth}
    \centering
    \includegraphics[width=.9\textwidth]{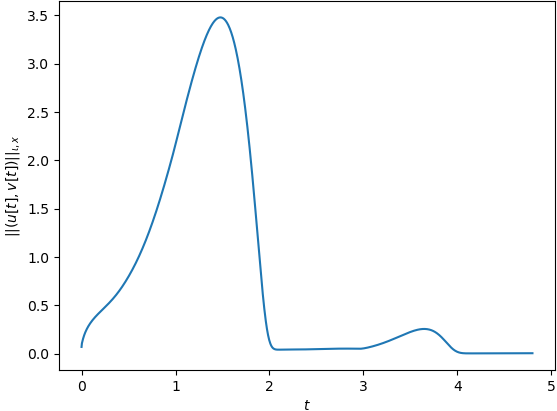}
    \caption{Closed loop. Control achieves stabilization. Propagation speed of 1 on the spatial interval $[0,1]$, results in an overall ``cycle'' of 2 seconds---it takes 2 seconds for the control input $U(\cdot)$ applied at $v( \cdot, 1)$ to reach $u( \cdot, 1)$.}
    \label{closed_loop_subfig}
\end{subfigure}
\begin{subfigure}[b]{\linewidth}
    \centering
    \includegraphics[width=0.9\linewidth]{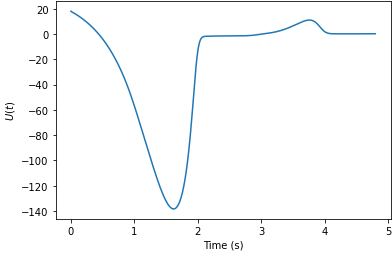} 
    \caption{
    Note the subintervals of action of control, $[0,2]$ and $[2,4]$, corresponding to the 2-second cycle of propagation from $x=1$ to 0 and back to 1. 
}
    \label{U_subfig}
\end{subfigure}
\caption{System \eqref{osysuiota-toy}--\eqref{osysbc-toy} in open loop (top) and under feedback (middle and bottom). 
See Figure \ref{visualisation} for how the control signal $U(t)$ achieves stabilization of the ensemble state $u(x,\iota,t)$.
}   
\end{figure}


\section{Numerical simulation}
\label{sec-simulations}

\begin{figure*}[p]
\centering
\includegraphics[width=0.43\linewidth]{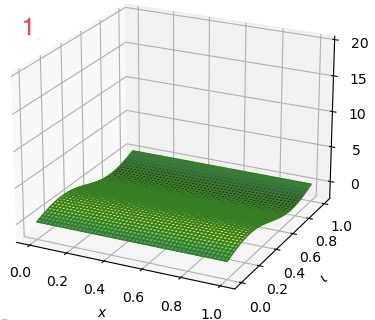}
\includegraphics[width=0.43\linewidth]{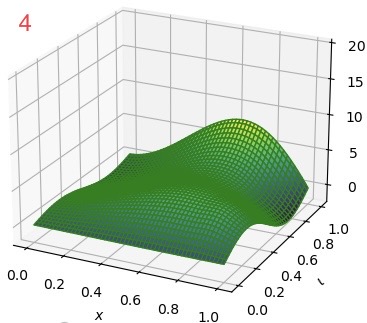}\\
\includegraphics[width=0.43\linewidth]{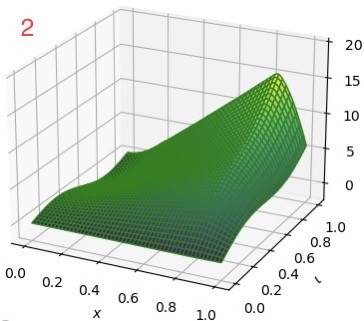}
\includegraphics[width=0.43\linewidth]{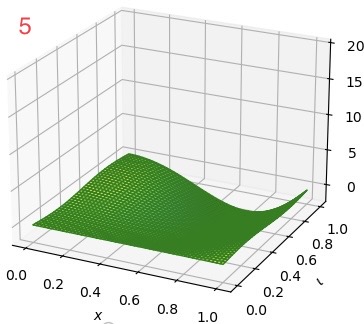}\\
\includegraphics[width=0.43\linewidth]{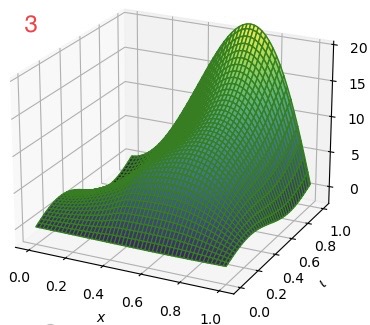}
\includegraphics[width=0.43\linewidth]{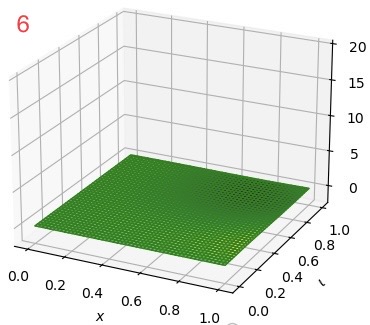}
\caption{Evolution of the unactuated {\em ensemble state} $u(t,x,\iota)$ in closed loop with scalar boundary control 
$v|_{x=1} =  \int_0^1 k(1,x)v(x)dx+  \int_0^1 \int_0^1 \kappa(1,x,\iota) u(\xi,\iota) d\iota dx$.
The time sequence is 1-2-3-4-5-6. The snapshorts are at times $t\in \{0, 0.16, 1.36, 1.82, 3.12, 6.62\}$. The system is open-loop unstable and one observes a large overshoot of the state at the time instant marked as 3. After that time, the control of the state $u$, though the boundary input on the state $v$, reaches the ensemble state $u$ and stabilizes it. The coupling in \eqref{osysuiota-toy} is most important at $\iota=1$ which is why the extremum of $u$ is always located along this line.}
\label{visualisation}
\end{figure*}

To illustrate the result of Theorem \ref{main_theorem} we numerically simulate the example system 
\begin{align}
 \label{osysuiota-toy}
u_t(t,x,y) + u_x(t,x,y)  
=& x(x+1) \left(y-\frac{1}{2}\right)\Bigg[e^xv(t,x) \nonumber\\ &
+x^2\int_0^1 \Big(\eta-\frac{1}{2}\Big) u(t,x,\eta)d\eta\Bigg]
 \nonumber \\&
\\
\label{osysvy-toy}
v_t(t,x) - v_x(t,x) 
=& -70e^{x\frac{35}{\pi^2}}\int_0^1y(y-1)u(t,x,y)dy
\\
\label{osysbc-toy}
u(t,0,y) =& \cos(2\pi y)v(t,0)\,
\end{align}
This system is a carefully constructed special version of \eqref{osysuiota}--\eqref{osysbound} for which the kernel PDE \eqref{keq1}--\eqref{kbound2}  is analytically solvable using a separation of variables and given by
\begin{eqnarray}
  k(x,\xi,y) &=& 35y(y-1)e^{2\xi\tilde{k}(x,\xi)}
  \\
  \tilde{k}(x,\xi)& =& \frac{35}{2\pi^2} 
\end{eqnarray}
   We choose such an example because the numerical solution of the kernel PDEs is a routine step, resolved by discretization in $y$ and the application of the techniques for a finite set of kernel PDEs already presented in \cite{dimeglio2013}. We choose an example in which we can highlight stabilization of a PDE ensemble, which is exponentially unstable in open loop, as seen in Fig \ref{open_loop_subfig}. In spite of the simplicity of the analytically solved kernel (independent of $x$ and quadratic in $y$), the feedback law (plotted in Figure \ref{U_subfig}) is not trivial and employs weighting of the state $u$ by the ensemble variable $y$:
\begin{align}
U(t) =& 35\int_0^1 \Bigg(e^{x{\frac{35}{\pi^2}}}\int_0^1 y(y-1)u(t,x,y)dy  
+\frac{1}{2\pi^2}v(t,x) \Bigg)dx
\end{align}
Equations \eqref{osysuiota-toy}--\eqref{osysbc-toy} are discretized in time and space using an explicit scheme where every function of $y$ are evaluated in 120 evenly spaced point.
The timestep is $\delta t = 0.004$ while the spatial domain is divided in 200 intervals of equal length. We run the closed loop system for 5s in order to observe the 
decay of the joint norm of the system which is indeed exponentially bounded (see Fig \ref{closed_loop_subfig}). This bound does not imply that $u$ is strictly decreasing as can be seen in Figure \ref{visualisation}.  

\section{Conclusions and Open Problems}
\label{sec-conclusions}

In this paper we presented initial results on stabilization, by boundary control, of \emph{ensemble} PDEs of the hyperbolic type. This concept is relevant for parametrized PDE collections with an inherently continuous parametrization. They differ from 2D PDEs in the sense that a derivative in the ensemble variable is absent from the model. The stabilization results are made possible by adapting the backstepping technique, a method that has already proven most effective for control of finite collections of hyperbolic PDEs. One can, in fact, infer the previous results for finite collections of hyperbolic PDEs by specializing the Hilbert space in the development we provide to a Euclidean space. We have had to develop a few  technical advances to enable the extension from the finite to infinite ensembles, including ensuring the regularity of the mathematical objects considered, especially the characteristic curves.

The problem of state estimation for finite ensembles of hyperbolic PDEs was tackled in \cite{dimeglio2013}, enabling also an observer-based output-feedback controller. Pursuing an \emph{ensemble} equivalent of such an observer result is a natural direction for further work on this topic. 

Of course, an extension of \cite{hu2016,hu2019boundary}, where $v$ in \eqref{osysviota} depends on a second ensemble variable, and the control, rather than being scalar-valued, is also a function of that additional ensemble variable, is also of interest.

The result of \cite{dimeglio2013} was extended in \cite[eqs. (1)-(5)]{deandrade2022backstepping} to coupled hyperbolic PDEs with a zero characteristic speed and an extension is of interest where the subsystem \cite[eqs. (3)]{deandrade2022backstepping} is an ensemble ODE.

For  ease of presentation, in this paper we consider the ensemble variable $y$ in the interval $[0,1]$. As long as $y$ stays within a compact there is no obvious technical reason why the results would not extend to multi-dimensional ensembles, but this is a topic for future research. We do note that, based on \cite{chen2022controllability}, controllability does not extend to higher-dimensional ensembles for ODEs. However, such a negative result for controllability  does not preclude stabilizability.

While in this paper we study stability in the $L_2$ sense, with respect to both the spatial variable $x$ and the ensemble variables $y$, it is of interest to explore stability in $L_\infty$ with respect to either or both of those variables. The second part of Theorem 2.3 in \cite{karafyllis2020stability} will be of significant help for that goal. 

Less immediate extensions would be in two additional directions: extensions to parabolic PDE ensembles and extensions to ensemble population dynamics---ensembles of nonlinear hyperbolic PDEs with in-domain actuation. 

Finally, we return to traffic as a motivating topic for an extension of PDE control frome finite ensembles \cite{burkhardt2021stop} to continuum ensembles. It is of interest to apply the present paper's design to such ARZ models of traffic with a continuum of classes. We point out that, in addition to the variability among cars and human driving styles, there will be variability in the headway settings in adaptive cruise control in automated vehicles, introducing one more ensemble variable and motivating the exploration for multivariable PDE ensembles. Furthermore, the two-dimensional Lighhill-Whitham-Richards (2D-LWR) model \cite{8569527} is a particular example of a hyperbolic ensemble PDE, with the ensemble variable being orthogonal to the net direction of traffic.


\appendix
\section{Appendix}

\subsection{Appendix: Stability of target system}
\label{target_system_stability}
\noindent\emph{Proof of Theorem \ref{targetstab}}:
Using proposition \ref{C0_kappa_C} we can consider 
\begin{equation}
M_{\boldsymbol{C}} = \sup_{(x,\xi) \in \mathcal{T}} \|\boldsymbol{C}(x,\xi)\|, \; M_{\kappa} = \sup_{(x,\xi) \in \mathcal{T}} \|\kappa[x,\xi]\|_E
\end{equation}  
Integrating by parts \eqref{Lyapounov} after using \eqref{osysu} and \eqref{osysv} we obtain
\begin{align}    
\dot{V}(t) =& \left[-pe^{-\delta x}||\alpha[t,x]||_E^2 + (1+x)\beta(t,x)^2 \right] _0^1 \nonumber
\\&- \int_0^1pe^{-\delta x} ||\alpha[t,x]||_E^2 + \beta(t,x)^2dx  \nonumber
\\ &+\int_0^1 2pe^{-\delta x} \bigg\langle \boldsymbol{L}(x)^{-1}\{\alpha[t,x]\}\;,\;  W[x]\beta(t,x) \nonumber \\ &+\int_0^x \boldsymbol{C}(x,\xi)\{\alpha[t,\xi]\} + \beta(t,\xi)\kappa[t,\xi] d\xi \nonumber
\\& +   \boldsymbol{\Theta}(x)\{\alpha[t,x]\} 
\bigg\rangle_E dx
\end{align}
 Using the inequality $ \langle a,b\rangle_E \le\frac{1}{2}  \left(||a||_E^2 + ||b||_E^2\right)$ we have
\begin{align}
&\int_0^1 2pe^{-\delta x} \Big\langle\boldsymbol{L}(x)^{-1}\{\alpha[t,x]\}\,,\, \boldsymbol{\Theta}(x)\{\alpha[t,x]\} + W[x]\beta(t,x) \nonumber
\\ &+ \int_0^x \boldsymbol{C}(x,\xi)\{\alpha(t,\xi)\} + \beta(t,\xi)\kappa[t,\xi] d\xi \Big\rangle_E dx \nonumber
\\ \nonumber
\le& \int_0^1 2pe^{-\delta x} \left\langle \boldsymbol{L}(x)^{-1}\circ \boldsymbol{\Theta}(x)\{\alpha[t,x]\} \,,\,\alpha[t,x] \right\rangle_E dx 
\\&+\int_0^1 pe^{-\delta x}\|W[x]\|_E^2\beta(t,x)^2 dx 
 \nonumber
\\ &+
\int_0^1 pe^{-\delta x} \left\|\int_0^x \boldsymbol{C}(x,\xi)\{\alpha(t,\xi)\}d\xi \right\|_E^2 dx \nonumber
\\&+ \int_0^1 pe^{-\delta x} \left\|
\int_0^x\beta(t,\xi)\kappa[t,\xi] d\xi \right\|_E^2dx \nonumber
\\ &+ \int_0^1 pe^{-\delta x}\left\| \boldsymbol{L}(x)^{-1}\{\alpha[t,x]\} \right\|_E^2dx
\end{align}
Therefore, if we introduce 
\begin{align}
    M_W = \sup_{x \in [0,1]} &\|W[x]\|_E,\; M_{\boldsymbol{L}^{-1}} = \sup_{x \in [0,1]} \|\boldsymbol{L}^{-1}(x)\| \\ \, &M_{\boldsymbol{\Theta}} = \sup_{x \in [0,1]} \|\boldsymbol{\Theta}(x)\|
\end{align}
we have
\begin{eqnarray}
 \dot{V}(t) &\le& (pe^{-\delta}\|q\|_E-1)\beta(t,0)^2 
 \nonumber \\
&&-\int_0^1\left(1-\frac{p}{\delta}M_{\kappa}^2e^{-\delta x} - pM_W^2e^{-\delta x}\right)\beta(t,x)^2dx
  \nonumber \\
&& -\int_0^1 pe^{-\delta x} \langle \alpha[t,x]\,,\, \boldsymbol{P}(x)\{\alpha[t,x]\} \rangle_E dx 
\end{eqnarray}
with \begin{equation}
    \boldsymbol{P} = \left(\delta - \frac{p}{\delta}M_{\boldsymbol{C}}^2\right)\boldsymbol{{\rm Id}} - 2 \boldsymbol{L}^{-1}\circ\boldsymbol{\Theta} - \boldsymbol{L}^{-1}
\end{equation}
Proceeding similarly as in \cite{dimeglio2013} we  choose $p$ small enough and $\delta$ big enough such that $\boldsymbol{P}$ is positive definite and the term multiplying $\beta(t,x)^2$ under the integral is positive. 
With $p$ no greater than one, and yet to be chosen,
we choose $\delta$ bigger than \begin{equation}
    \delta^* =1+  M_{\boldsymbol{C}}^2 + 2M_{\boldsymbol{L}^{-1}}M_{\boldsymbol{\Theta}} + M_{\boldsymbol{L}^{-1}}\,\end{equation} for
the operator $\boldsymbol{P}$ to be positive definite. This is seen by computing
\begin{align}
    \forall x \in [0,1],\,\forall a \in E, \,\nonumber\\
    \langle a\,,\, P(x)\{a\}\rangle_E \ge& \left\langle a\,,\, \delta a\right\rangle_E - \big\langle a\,,\, (M_{\boldsymbol{C}}^2\boldsymbol{\rm {Id}} \nonumber \\ &+ 2 \boldsymbol{L}(x)^{-1}\circ\boldsymbol{\Theta}(x) + \boldsymbol{L}(x)^{-1})\ \big\rangle_E
    \\
    \ge& (\delta - \delta^*) \|a\|_E^2
\end{align}
Once $\delta > \delta^*$ is chosen, we select $p$ such that
\begin{equation}
    0 < p < \min \left(1, \frac{e^{\delta}}{\|q\|_E}, \frac{\delta}{M_{\kappa}^2+\delta M_W^2}\right)
\end{equation}
which ensures that
\begin{eqnarray}
&&\int_0^1\left(1-\frac{p}{\delta}M_{\kappa}^2e^{-\delta x} - pM_W^2e^{-\delta x}\right)\beta(t,x)^2dx \nonumber \\ &\ge& \int_0^1\left(1-p\left(\frac{M_{\kappa}^2}{\delta} + M_W^2\right)\right)\beta(t,x)^2dx\\
  &&(pe^{-\delta}\|q\|_E-1)\beta(t,0)^2 \le 0\,
\end{eqnarray}
and hence, 
for any $\varepsilon$ verifying 
\begin{equation}
    0<\varepsilon \le \min\left(p(\delta - \delta^*), \frac{p}{\delta}(M_{\kappa}^2+\delta M_W^2)\right)
\end{equation}
we have that
\begin{eqnarray}
\label{ineqnorme}
   \dot{V}(t) \le - \varepsilon \left\| \left( \begin{array}{c}
      \alpha[t]\\
      \beta[t]
\end{array} \right) \right\|_{x,y}^2 \,,
\end{eqnarray}
or, in the scalar-valued notation, 
\begin{equation}
    \dot V(t) \leq - \varepsilon\int_0^1 \left[\int_0^1 \alpha^2(t,x,y) dy + \beta^2(t,x)\right] dx. 
\end{equation}

We then need to show that $\sqrt{V}$ is equivalent to $\|\cdot\|_{x,y}$ in order to conclude the proof. If we consider $m_{\mu} = \min_x(\mu)$ and $M_{\mu} = \max_x(\mu)$ we have
\begin{eqnarray}
\label{equivalent_norm_left}
 \min\left(\frac{pe^{-\delta}}{M_{\boldsymbol{L}^{-1}}}, \, \frac{1}{m_{\mu}}\right) \left\| \left( \begin{array}{c}
      \alpha[t]\\
      \beta[t]
\end{array} \right) \right\|_{x,y}^2 \le  V(t) \\\label{equivalent_norm_right}
V(t) \le 
\left( p M_{\boldsymbol{L}^{-1}}+\frac{2}{M_{\mu}} \right) \left\| \left( \begin{array}{c}
      \alpha[t]\\
      \beta[t]
\end{array} \right) \right\|_{x,y}^2
\end{eqnarray}
Combining \eqref{equivalent_norm_left}-\eqref{equivalent_norm_right} and \eqref{ineqnorme} we get
\begin{align}
& \left\| \left( \begin{array}{c}
      \alpha[t]\\
      \beta[t]
\end{array} \right) \right\|_{x,y} \le M_1
\left\| \left( \begin{array}{c}
      \alpha[0]\\
      \beta[0]
\end{array} \right) \right\|_{x,y} e^{-\tau t} \\
\text{where}\;& \qquad M_1 = \sqrt{ \frac{p M_{\boldsymbol{L}^{-1}}+\frac{2}{M_{\mu}}}{ \min\left(\frac{pe^{-\delta}}{M_{\boldsymbol{L}^{-1}}}, \, \frac{1}{m_{\mu}}\right)}} \\
& \qquad  \tau = \frac{\varepsilon 
}{2\min\left(\frac{pe^{-\delta}}{M_{\boldsymbol{L}^{-1}}}, \, \frac{1}{m_{\mu}}\right)} \end{align}
This completes the proof. 
\hfill $\blacksquare$


\subsection{Appendix : Derivation of kernel equation
} \label{Appendix_kernel_design}
We have
\begin{align}
\beta_x(t,x) =& v_x(t,x) - \langle k[x,x]\,,\, u[t,x]\rangle_E - \Tilde{k}(x,x)v(t,x) \nonumber
\\ & - \int_0^x \langle k_x[x,\xi]\,,\, u[t,\xi]\rangle_E d\xi - \int_0^x \Tilde{k}_x(x,\xi)v(t,\xi)d\xi
\end{align}
and using \eqref{osysu}-\eqref{osysv}
\begin{align}
\beta_t(t,x) =& v_t(t,x) - \int_0^x \langle k[x,\xi]\,,\, u_t[t,\xi]\rangle_E d\xi \nonumber\\ &- \int_0^x \Tilde{k}(x,\xi)v_t(t,\xi)d\xi\\
=& \mu(x)v_x(t,x) + \langle \Xi[x]\; , \; u[t,x] \rangle_E \nonumber\\ &- \int_0^x \Big\langle k[x,\xi]\,,\, -\boldsymbol{L}(\xi)\{u_x[t,\xi]\} + \boldsymbol{\Theta}(\xi)\{u[t,\xi]\} \nonumber \\&+ W[\xi]v[t,\xi]\Big\rangle_E d\xi \nonumber \\ &-\int_0^x \Tilde{k}(x,\xi)\left(\mu(\xi)v_x(t,\xi) + \left\langle \Xi[\xi]\; , \; u[t,\xi] \right \rangle_E\right)d\xi
\end{align}
Plugging both expressions in \eqref{betasys} and integrating by part we get 
 \begin{align}
 0 =& \mu(x)v_x(t,x) + \langle \Xi[x]\, ,\, u[t,x] \rangle_E 
 \nonumber -\int_0^x \bigg \langle \boldsymbol{L}'(\xi)\{k[x,\xi]\} \nonumber \\ &+ \boldsymbol{L}(\xi)\{k_\xi[x,\xi]\} + \boldsymbol{\Theta}^t(\xi)\{k[x,\xi]\} +  \Tilde{k}(x,\xi)\Xi[\xi] \nonumber
 \\ &  - \mu(x)k_x[x,\xi]\,,\, u[t,\xi] \bigg \rangle_E d\xi- \left[ \Tilde{k}(x,\xi)\mu(\xi)v(t,\xi)\right]_0^x 
 \nonumber \\
 &- \int_0^x \bigg ( \langle W[\xi],k[x,\xi]  \rangle_E - \Tilde{k}_\xi(x,\xi)\mu(\xi)-\Tilde{k}(x,\xi)\mu'(\xi) 
  \nonumber 
 \\&- \Tilde{k}_x(x,\xi)\mu(x) \bigg ) v(t,\xi) d\xi
 -\mu(x)v_x(t,x)
\nonumber \\
& + \mu(x) \langle k[x,x], u[t,x] \rangle_E + \mu(x) \Tilde{k}(x,x)v(t,x)
\nonumber \\
&+ \left[ \Big \langle \boldsymbol{L}(\xi)\{k[x,\xi]\}\,,\, u[t,\xi] \Big\rangle_E \right]_0^x
\\
=& \langle \Xi[x] \,,\, u[t,x] \rangle_E
\nonumber 
- \int_0^x \bigg \langle \boldsymbol{L}'(\xi)\{k[x,\xi]\}\nonumber \\&+ \boldsymbol{L}(\xi)\{k_\xi[x,\xi]\} 
 + \boldsymbol{\Theta}^t(\xi)\{k[x,\xi]\} +\Tilde{k}(x,\xi)\Xi[\xi] \nonumber\\&   - \mu(x)k_x[x,\xi]\,,\, u[t,\xi] \bigg \rangle_E d\xi + \Tilde{k}(x,0)\mu(0)v(t,0)
\nonumber\\&  - \int_0^x \Big ( \langle W[\xi]\,,\,k[x,\xi]  \rangle_E
- \Tilde{k}_\xi(x,\xi)\mu(\xi)-\Tilde{k}(x,\xi)\mu'(\xi) 
\nonumber \\
& - \Tilde{k}_x(x,\xi)\mu(x) \Big ) v(t,\xi) d\xi+\mu(x) \langle k[x,x]\,,\, u[t,x] \rangle_E 
\nonumber\\
& +\left[ \Big \langle \boldsymbol{L}(\xi)\{k[x,\xi]\}\,,\, u[t,\xi] \Big\rangle_E \right]_0^x \label{end_equation}
 \end{align}
One can then see that a solution $(k,\Tilde{k})$ of the system \eqref{keq1}-\eqref{kbound2} would verify \eqref{end_equation}.
\subsection{Appendix: Inverse kernel}
\label{inverse_kernel_section}We introduce
\begin{eqnarray}
  \Gamma(t,x) = \beta(t,x) + \int_0^x \langle k[x,\xi]\,,\, \alpha[t,\xi]\rangle_E d\xi
\end{eqnarray} in order to write
\begin{eqnarray}
   \label{volterra_v}
   v(t,x) &=& \int_0^x \Tilde{k}(x,\xi)v(t,x,\xi)d\xi +\Gamma(t,x) 
\end{eqnarray}
Notice that \eqref{volterra_v} can be seen as a scalar Volterra equation of the second kind, similar to \eqref{Volterra_kappa}, with scalar unknown $v$ instead of vector $\kappa$.
Using \cite[ch. 8]{vrabieDifferentialEquationsIntroduction2016} / adapting Lemma \ref{volterra_solution} we have the existence of the continuous inverse kernel $\Tilde{l}$ (similar to $\tilde{k}^{\infty}$ in Lemma \ref{volterra_solution}) such that
\begin{align}
v(t,x) =& \Gamma(t,x) + \int_0^x \Tilde{l}(x,\xi)\Gamma(t,\xi)d\xi\\
=&\Gamma(t,x) + \int_0^x \Tilde{l}(x,\xi) \Bigg(  \beta(t,\xi) \nonumber \\ &+ \int_0^\xi \langle k[\xi,s], \alpha[t,s]\rangle_E ds\Bigg)d\xi
\\
=&\Gamma(t,x) + \int_0^x \Tilde{l}(x,\xi) \beta(t,\xi)d\xi \nonumber\\&+ \int_0^x \int_{s}^x \Tilde{l}(x,\xi) \langle k[\xi,s], \, \alpha[t,s]\rangle_E d\xi ds
\\
=&\Gamma(t,x) \nonumber
+ \int_0^x \Bigg( \Tilde{l}(x,\xi) \beta(t,\xi) \\&+ \left\langle \int_{\xi}^x \Tilde{l}(x,s)  k[s, \xi] ds, \, \alpha[t,\xi] \right \rangle_E \Bigg)  d\xi
\\
=&\beta(t,x) + \int_0^x \Tilde{l}(x,\xi) \beta(t,\xi)d\xi \nonumber \\ &+ \int_0^x \langle l[x,\xi], \, \alpha[t,\xi] \rangle_E d\xi 
\\ \text{where} \qquad &l[x,\xi] = k[x,\xi] + \int_{\xi}^x \Tilde{l}(x,s)  k[s, \xi] ds
\end{align}




\section*{Acknowledgments}

We thank Jean Auriol, Xudong Chen, Florent Di Meglio, Rafael Vazquez, and Huan Yu for valuable comments on a draft of this paper. 

This work was supported in part by the Air Force Office of Scientific Research under grant FA9550-22-1-0265. 

\bibliographystyle{IEEEtranS}
\bstctlcite{IEEEexample:BSTcontrol}
\bibliography{refs-Krstic}

\begin{thebibliography}{10}
\providecommand{\url}[1]{#1}
\csname url@samestyle\endcsname
\providecommand{\newblock}{\relax}
\providecommand{\bibinfo}[2]{#2}
\providecommand{\BIBentrySTDinterwordspacing}{\spaceskip=0pt\relax}
\providecommand{\BIBentryALTinterwordstretchfactor}{4}
\providecommand{\BIBentryALTinterwordspacing}{\spaceskip=\fontdimen2\font plus
\BIBentryALTinterwordstretchfactor\fontdimen3\font minus
  \fontdimen4\font\relax}
\providecommand{\BIBforeignlanguage}[2]{{%
\expandafter\ifx\csname l@#1\endcsname\relax
\typeout{** WARNING: IEEEtranS.bst: No hyphenation pattern has been}%
\typeout{** loaded for the language `#1'. Using the pattern for}%
\typeout{** the default language instead.}%
\else
\language=\csname l@#1\endcsname
\fi
#2}}
\providecommand{\BIBdecl}{\relax}
\BIBdecl

\bibitem{Anfinsen2019Adaptive}
H.~Anfinsen and O.~Aamo, \emph{Adaptive Control of Hyperbolic {PDE}s}.\hskip
  1em plus 0.5em minus 0.4em\relax Springer, 2019.

\bibitem{Auriol2018Delay}
J.~Auriol, U.~Aarsnes, P.~Martin, and F.~Di{ }Meglio, ``Delay-robust control
  design for two heterodirectional linear coupled hyperbolic {PDE}s,''
  \emph{IEEE Transactions on Automatic Control}, vol.~63, no.~10, pp.
  3551--3557, 2018.

\bibitem{Auriol2018Delay1}
J.~Auriol, F.~Bribiesca-Argomedo, D.~Saba, M.~D. Loreto, and F.~Di{ }Meglio,
  ``Delay-robust stabilization of a hyperbolic {PDE-ODE} system,''
  \emph{Automatica}, vol.~95, pp. 494--502, 2018.

\bibitem{BERNARD20142692}
P.~Bernard and M.~Krstic, ``Adaptive output-feedback stabilization of non-local
  hyperbolic pdes,'' \emph{Automatica}, vol.~50, no.~10, pp. 2692--2699, 2014.

\bibitem{burkhardt2021stop}
M.~Burkhardt, H.~Yu, and M.~Krstic, ``Stop-and-go suppression in two-class
  congested traffic,'' \emph{Automatica}, vol. 125, p. 109381, 2021.

\bibitem{chen2022backstepping}
G.~Chen, R.~Vazquez, and M.~Krstic, ``Backstepping-based exponential
  stabilization of timoshenko beam with prescribed decay rate,''
  \emph{IFAC-PapersOnLine}, vol.~55, no.~26, pp. 162--167, 2022.

\bibitem{chen2019structure}
X.~Chen, ``Structure theory for ensemble controllability, observability, and
  duality,'' \emph{Mathematics of Control, Signals, and Systems}, vol.~31,
  no.~2, pp. 1--40, 2019.

\bibitem{chen2021sparse}
X.~Chen, ``Sparse linear ensemble systems and structural controllability,''
  \emph{IEEE Transactions on Automatic Control}, vol.~67, no.~7, pp.
  3337--3348, 2021.

\bibitem{Coron2013Local}
J.~Coron, R.~Vazquez, M.~Krstic, and G.~Bastin, ``Local exponential {$H^2$}
  stabilization of a $2\times2$ quasilinear hyperbolic system using
  backstepping,'' \emph{SIAM Journal on Control and Optimization}, vol.~51,
  no.~3, pp. 2005--2035, 2013.

\bibitem{deandrade2022backstepping}
\BIBentryALTinterwordspacing
G.~A. de~Andrade, R.~Vazquez, I.~Karafyllis, and M.~Krstic, ``Backstepping
  control of a hyperbolic pde system with zero characteristic speed states,''
  2022. [Online]. Available: \url{https://arxiv.org/abs/2211.14290}
\BIBentrySTDinterwordspacing

\bibitem{deutscher2018}
J.~Deutscher and J.~Gabriel, ``Minimum time output regulation for general
  linear heterodirectional hyperbolic systems,'' \emph{International Journal of
  Control}, vol.~93, pp. 1826--1838, 2018.

\bibitem{DEUTSCHER201556}
J.~Deutscher, ``A backstepping approach to the output regulation of boundary
  controlled parabolic pdes,'' \emph{Automatica}, vol.~57, pp. 56--64, 2015.

\bibitem{dimeglio2013}
F.~Di~Meglio, R.~Vazquez, and M.~Krstic, ``Stabilization of a system of $n+1$
  coupled first-order hyperbolic linear {PDE}s with a single boundary input,''
  \emph{IEEE Transactions on Automatic Control}, vol.~58, pp. 3097--3111, 2013.

\bibitem{diagne2017backstepping}
A.~Diagne, M.~Diagne, S.~Tang, and M.~Krstic, ``Backstepping stabilization of
  the linearized saint-venant--exner model,'' \emph{Automatica}, vol.~76, pp.
  345--354, 2017.

\bibitem{diagne2017control}
M.~Diagne, S.-X. Tang, A.~Diagne, and M.~Krstic, ``Control of shallow waves of
  two unmixed fluids by backstepping,'' \emph{Annual Reviews in Control},
  vol.~44, pp. 211--225, 2017.

\bibitem{dirr2021uniform}
G.~Dirr and M.~Sch{\"o}nlein, ``Uniform and lq-ensemble reachability of
  parameter-dependent linear systems,'' \emph{Journal of Differential
  Equations}, vol. 283, pp. 216--262, 2021.

\bibitem{DIMEGLIO2018281}
F.~Di Meglio, F.~B. Argomedo, L.~Hu, and M.~Krstic, ``Stabilization of coupled
  linear heterodirectional hyperbolic pde–ode systems,'' \emph{Automatica},
  vol.~87, pp. 281--289, 2018.

\bibitem{chen2023backstepping}
C.~Guangwei, V.~Rafael, and K.~Miroslav, ``Backstepping-based rapid
  stabilization of two-layer timoshenko composite beams,'' \emph{IFAC World
  Congress}, 2022.

\bibitem{helmke2014uniform}
U.~Helmke and M.~Sch{\"o}nlein, ``Uniform ensemble controllability for
  one-parameter families of time-invariant linear systems,'' \emph{Systems \&
  Control Letters}, vol.~71, pp. 69--77, 2014.

\bibitem{hu2016}
L.~Hu, F.~Di~Meglio, R.~Vazquez, and K.~M., ``Control of homodirectional and
  general heterodirectional linear coupled hyperbolic {PDE}s,'' \emph{IEEE
  Transactions on Automatic Control}, vol.~61, pp. 3301--3314, 2016.

\bibitem{hu2019boundary}
L.~Hu, R.~Vazquez, F.~Di{ }Meglio, and M.~Krstic, ``Boundary exponential
  stabilization of 1-dimensional inhomogeneous quasi-linear hyperbolic
  systems,'' \emph{SIAM J. Control and Optimization}, vol.~57, no.~2, pp.
  963--998, 2019.

\bibitem{Yu2022}
Y.~Huan and K.~Miroslav, \emph{Traffic Congestion Control by PDE
  Backstepping}.\hskip 1em plus 0.5em minus 0.4em\relax Springer, 2022.

\bibitem{karafyllis2020stability}
I.~Karafyllis and M.~Krstic, ``Stability results for the continuity equation,''
  \emph{Systems \& Control Letters}, vol. 135, p. 104594, 2020.

\bibitem{krstic2008Backstepping}
M.~Krstic and A.~Smyshlyaev, ``Backstepping boundary control for first-order
  hyperbolic {PDE}s and application to systems with actuator and sensor
  delays,'' \emph{Systems $\&$ Control Letters}, vol.~57, no.~9, pp. 750--758,
  2008.

\bibitem{li2010ensemble}
J.-S. Li, ``Ensemble control of finite-dimensional time-varying linear
  systems,'' \emph{IEEE Transactions on Automatic Control}, vol.~56, no.~2, pp.
  345--357, 2010.

\bibitem{li2009ensemble}
J.-S. Li and N.~Khaneja, ``Ensemble control of bloch equations,'' \emph{IEEE
  Transactions on Automatic Control}, vol.~54, no.~3, pp. 528--536, 2009.

\bibitem{8569527}
S.~Mollier, M.~L. Delle~Monache, and C.~Canudas-de Wit, ``2d-lwr in large-scale
  network with space dependent fundamental diagram,'' in \emph{2018 21st
  International Conference on Intelligent Transportation Systems (ITSC)}, 2018,
  pp. 1640--1645.

\bibitem{qi2019stabilization}
J.~Qi, M.~Krstic, and S.~Wang, ``Stabilization of reaction--diffusions pde with
  delayed distributed actuation,'' \emph{Systems \& Control Letters}, vol. 133,
  p. 104558, 2019.

\bibitem{vazquez2016boundary}
V.~Rafael and K.~Miroslav, ``Boundary control of coupled
  reaction-advection-diffusion systems with spatially-varying coefficients,''
  \emph{IEEE Transactions on Automatic Control}, vol.~62, no.~4, pp.
  2026--2033, 2016.

\bibitem{fedaabdelazizmustafaAnalyticalNumericalSolutions2014}
F.~A. A.~M. Salameh, ``Analytical and {{Numerical Solutions}} of {{Volterra
  Integral Equation}} of the {{Second Kind}},'' Ph.D. dissertation, An-Najah
  National University, 2014.

\bibitem{1369395}
A.~Smyshlyaev and M.~Krstic, ``Closed-form boundary state feedbacks for a class
  of 1-d partial integro-differential equations,'' \emph{IEEE Transactions on
  Automatic Control}, vol.~49, no.~12, pp. 2185--2202, 2004.

\bibitem{vazquez2006explicit}
R.~Vazquez and M.~Krstic, ``Explicit integral operator feedback for local
  stabilization of nonlinear thermal convection loop pdes,'' \emph{Systems \&
  control letters}, vol.~55, no.~8, pp. 624--632, 2006.

\bibitem{vazquez2007closed}
R.~Vazquez and M.~Krstic, ``A closed-form feedback controller for stabilization
  of the linearized 2-d navier--stokes {P}oiseuille system,'' \emph{IEEE
  Transactions on Automatic Control}, vol.~52, no.~12, pp. 2298--2312, 2007.

\bibitem{vrabieDifferentialEquationsIntroduction2016}
I.~I. Vrabie, \emph{Differential Equations: An Introduction to Basic Concepts,
  Results and Applications}, 3rd~ed.\hskip 1em plus 0.5em minus 0.4em\relax
  {World Scientific}, 2016.

\bibitem{WANG2020109131}
J.~Wang and M.~Krstic, ``Delay-compensated control of sandwiched
  ode–pde–ode hyperbolic systems for oil drilling and disaster relief,''
  \emph{Automatica}, vol. 120, p. 109131, 2020.

\bibitem{9319184}
J.~Wang and M.~Krstic, ``Event-triggered output-feedback backstepping control
  of sandwich hyperbolic pde systems,'' \emph{IEEE Transactions on Automatic
  Control}, vol.~67, no.~1, pp. 220--235, 2022.

\bibitem{xu2008stabilization}
C.~Xu, E.~Schuster, R.~Vazquez, and M.~Krstic, ``Stabilization of linearized 2d
  magnetohydrodynamic channel flow by backstepping boundary control,''
  \emph{Systems \& control letters}, vol.~57, no.~10, pp. 805--812, 2008.

\bibitem{chen2020ensemble}
C.~Xudong, ``Ensemble observability of bloch equations with unknown population
  density,'' \emph{Automatica}, vol. 119, p. 109057, 2020.

\bibitem{chen2022controllability}
\BIBentryALTinterwordspacing
C.~Xudong, ``Controllability issues of linear ensemble systems over
  multi-dimensional parameterization spaces,'' 2022. [Online]. Available:
  \url{https://arxiv.org/abs/2003.04529}
\BIBentrySTDinterwordspacing

\bibitem{Yu2019Traffic}
H.~Yu and M.~Krstic, ``Traffic congestion control for {Aw-Rascle-Zhang}
  model,'' \emph{Automatica}, vol. 100, pp. 38--51, 2019.

\end{thebibliography}

 \newpage
 \begin{IEEEbiography}[{\includegraphics[width=1in,height=1.25in,clip,keepaspectratio]{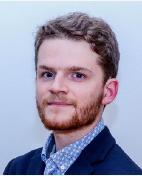}}]
{Valentin Alleaume}
is a graduate student at Ecole des Mines de Paris, also known as Mines Paris-PSL, France, graduating with a Master’s degre in 2024. He was a visiting student at University of California, San Diego, in summer of 2022.

\end{IEEEbiography} 

\newpage
\begin{IEEEbiography}
[{\includegraphics[width=1in,height=1.25in,clip,keepaspectratio]{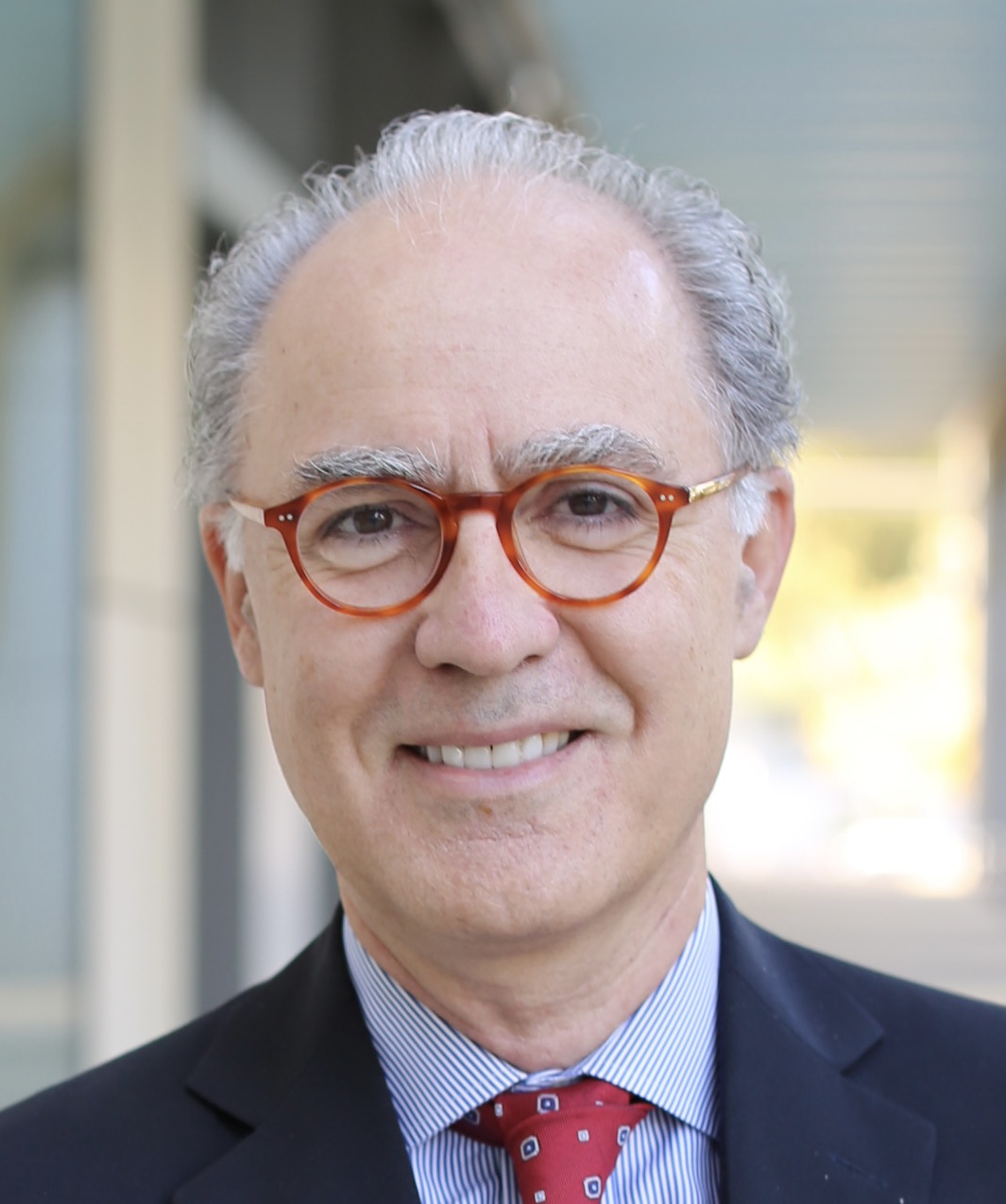}}]
{Miroslav Krstic}
 is distinguished professor, Alspach endowed chair, founding director of the  Center for Control Systems and Dynamics, and senior associate vice chancellor for research at UC San Diego. 
 He is Fellow of IEEE, IFAC, ASME, SIAM, AAAS, IET (UK), and AIAA (Assoc. Fellow), as well as foreign member of the Serbian Academy of Sciences and Arts and of the Academy of Engineering of Serbia. He has received the Bellman Award, Bode Lecture Prize, SIAM Reid Prize, ASME Oldenburger Medal, Nyquist Lecture Prize, Paynter Award, Ragazzini  Award, IFAC Nonlinear Control Systems Award, Chestnut Award, IFAC Distributed Parameter Systems Award, IFAC Adaptive and Learning Systems Award, AV Balakrishnan Award for the Mathematics of Systems, CSS Distinguished Member Award, the PECASE, NSF Career, and ONR YI awards, the Schuck (1996 and 2019) and Axelby award, and the first UCSD Research Award given to an engineer. 
 He is EiC of Systems \& Control Letters and has been senior editor in Automatica, IEEE Trans. Automatic Control, and of two Springer book series. 
 Krstic has coauthored eighteen books on adaptive, nonlinear, and stochastic control, extremum seeking, control of PDE systems including turbulent flows, and control of delay systems.
\end{IEEEbiography} 

\vfill

\end{document}